\documentclass[envcountsame,orivec,runningheads,11pt,a4paper]{llncs}
\usepackage{graphics,xspace,amsmath,amssymb,amsfonts,xspace,stmaryrd,textgreek,graphicx}
\usepackage[utf]{kotex}
\usepackage{hyperref,breakurl}             
\usepackage[normalem]{ulem}
\usepackage[capitalise]{cleveref}
\usepackage{a4}
\usepackage{tensor}
\newcommand{\IR}{\mathbb{R}}

\newcommand{\IC}{\mathbb{C}}

\newcommand{\IN}{\mathbb{N}}
\newcommand{\poly}{\operatorname{poly}}

\newcommand{\LIP}{\operatorname{Lip}_1}
\newcommand{\reduceq}{\preccurlyeq}
\newcommand{\reduceqL}{\preccurlyeq_{\rm L}}
\newcommand{\reduceqP}{\preccurlyeq_{\rm P}}
\newcommand{\Card}{\operatorname{Card}}
\newcommand{\bin}{\mathrm{bin}}
\newcommand{\sbin}{\mathrm{sbin}}

\newcommand{\dom}{\operatorname{dom}}
\newcommand{\range}{\operatorname{range}}

\newcommand{\id}{\operatorname{id}}

\newcommand{\sdzero}{\textup{\texttt{0}}\xspace}
\newcommand{\sdone}{\textup{\texttt{1}}\xspace}

\newcommand{\calB}{\mathcal{B}}
\newcommand{\calC}{\mathcal{C}}

\newcommand{\calH}{\mathcal{H}}

\newcommand{\calO}{\mathcal{O}}
\newcommand{\calP}{\mathcal{P}}
\newcommand{\calK}{\mathcal{K}}
\newcommand{\calU}{\mathcal{U}}
\newcommand{\calV}{\mathcal{V}}
\newcommand{\calY}{\mathcal{Y}}

\newcommand{\Hilbert}{\mathcal{H}}
\newcommand{\TWO}{\ensuremath{\{\sdzero,\sdone\}}}
\newcommand{\THREE}{\ensuremath{\{-1,0,1\}}}
\newcommand{\unary}[1]{\widehat{#1}}
\newcommand{\dual}[1]{#1^\dagger}   
\newcommand{\Cantor}{\calC}
\newcommand{\Baire}{\calB}
\newcommand{\tildedelta}{\widetilde{\delta}}
\newcommand{\tildesbin}{\widetilde{\operatorname{sbin}}}
\newcommand{\deltabox}{\delta _{\square}}

\newcommand{\toto}{\rightrightarrows}
\newcommand{\mapstoto}{\Mapsto}
\newcommand{\loinv}[1]{#1^{\underline{-1}}}
\newcommand{\App}{\operatorname{App}}
\newcommand{\diam}{\operatorname{diam}}

\newcommand{\cball}{\overline{\operatorname{B}}}

\newcommand{\dfeq}{:=}
\newcommand{\pcolon}{\mathpunct{\,:\subseteq}}

\newcommand{\ie}{i.\,e.\xspace}

\newtheorem{observation}[theorem]{Observation}

\newtheorem{fact}[theorem]{Fact}
\newcommand{\Entropy}{\eta}
\newcommand{\Capacity}{\kappa}
\newcommand{\ENTROPY}{\calH}
\newcommand{\CAPACITY}{\calK}
\newcommand{\myoplus}[2]{\!\tensor[_{#1}]{\oplus\!\!}{_{#2}}\:}
\newcommand{\myOplus}[2]{\bigoplus_{#2}^{#1}}
\newcommand{\arXiv}[1]{\href{http://arXiv.org/abs/#1}{\texttt{[arXiv:#1]}}\xspace}
\newcommand{\mydoi}[1]{\href{http://doi.org/#1}{\texttt{[doi:#1]}\xspace}}
\newcommand{\COMMENTED}[1]{}
\usepackage[usenames]{color}  
\newcommand{\Martin}[1]{#1}

\makeatletter
\newcommand*\bigdot{\mathpalette\bigdot@{.5}}
\newcommand*\bigdot@[2]{\mathbin{\vcenter{\hbox{\scalebox{#2}{$\m@th#1\bullet$\,\,}}}}}
\makeatother

\title{Quantitative Coding and Complexity Theory \\ of \emph{Continuous} Data\thanks{%
Supported by the National Research Foundation of Korea 
(grant NRF-2017R1E1A1A03071032)
and the International Research \& Development Program of
the Korean Ministry of Science and ICT (grant NRF-2016K1A3A7A03950702).
This work has grown from preprints \arXiv{2002.04005} and \arXiv{1809.08695}, from \cite{Donghyun}, 
and from \cite{DBLP:conf/cie/Lim020}: based on discussions 
with Bruce Kapron, Akitoshi Kawamura, Sewon Park, Gleb Pogudin, Matthias Schr\"{o}der, and Florian Steinberg.}}
\titlerunning{Quantitative Coding and Complexity Theory of Metric Spaces}
\institute{KAIST, School of Computing, Republic of Korea}
\author{Donghyun Lim \and Martin Ziegler}
\date{}

\begin{document}
\maketitle

\begin{abstract}
Specifying a computational problem includes fixing encodings for input and output:
encoding graphs as adjacency matrices, characters as integers, integers as bit strings, and vice versa.
For such discrete data, the actual encoding is usually straightforward and/or complexity-theoretically inessential (up to linear or polynomial time, say).
Concerning continuous data, already real numbers naturally suggest various encodings 
(formalized as historically so-called representations) with very different properties,
ranging from the computably `unreasonable' binary expansion \mydoi{10.1112/plms/s2-43.6.544} via qualitatively to polynomially and even linearly complexity-theoretically `reasonable' signed-digit expansion. 
But how to distinguish between un/suitable encodings of other spaces common in Calculus and Numerics, such as Sobolev?

With respect to qualitative computability over topological spaces, \emph{admissibility} had been identified \mydoi{10.1016/0304-3975(85)90208-7} 
as a crucial criterion for a representation over the Cantor space of infinite binary sequences to be `reasonable':
It requires the representation to be (sequentially) continuous, and to be maximal with respect to (sequentially) continuous reduction \mydoi{10.1007/11780342\_48}. 
Such representations are guaranteed to exist for a large class of spaces.
And for (precisely) these does the sometimes so-called \emph{Main Theorem} hold: 
which characterizes continuity of functions by the continuity of mappings translating codes, so-called \emph{realizers}. 

We refine qualitative computability over topological spaces to quantitative complexity over metric spaces,
by developing the theory of \emph{polynomially} and of \emph{linearly admissible} representations.
Informally speaking, these are `optimally' continuous, namely linearly/polynomially relative to the space's entropy; and maximal with respect to relative linearly/polynomially continuous reductions defined below.
A large class of spaces is shown to admit a quantitatively admissible representation, including a generalization of the signed-digit encoding;
and these exhibit a quantitative strengthening of the qualitative \emph{Main Theorem},
namely now characterizing quantitative continuity of functions by quantitative continuity of realizers.
Our quantitative admissibility thus provides the desired criterion for complexity-theoretically `reasonable' encodings.

We then rephrase quantitative admissibility as quantitative continuity of both the representation and of its set-valued inverse,
the latter adopting from \mydoi{10.4115/jla.2013.5.7} 
a new notion of `sequential' continuity for multifunctions.
By establishing a quantitative continuous selection theorem for multifunctions between compact ultrametric spaces, we can
extend our above quantitative \emph{Main Theorem} from functions to multifunctions aka search problems.
Higher-type complexity is captured by generalizing Cantor's (and Baire's) ground space for encodings to other (compact) \emph{ultra}metric spaces.
\end{abstract}

\setcounter{tocdepth}{3}
\renewcommand{\contentsname}{}
\begin{center}
\begin{minipage}[c]{0.98\textwidth}\vspace*{-9ex}%
\tableofcontents
\end{minipage}
\end{center}

\section{Introduction}

Machine models formalize computation:
Such a model specifies the ways to input, operate on, and output elements from some fixed `ground' space $\calU$;
as well as a measure of cost, and of parameter(s) in whose dependence to analyze and bound said computational cost.
Computational problems over spaces $X$ other than $\calU$ are treated by encoding its elements/instances over $\calU$,
formalized as \emph{representations}: surjective partial mappings $\xi:\subseteq\calU\twoheadrightarrow X$.

Machine models formalize computation:
Such a model specifies the ways to input, operate on, and output elements from some fixed `ground' space $\calU$;
as well as a measure of cost, and of parameter(s) in whose dependence to analyze and bound said computational cost.
Computational problems over spaces $X$ other than $\calU$ are treated by encoding its elements/instances over $\calU$,
formalized as \emph{representations}: surjective partial mappings $\xi:\subseteq\calU\twoheadrightarrow X$.

\begin{example}
\label{x:Example}
\begin{enumerate}
\item[a)]
Recall the Turing machine model operating on the set $\calU=\TWO^*$ of finite binary sequences.
Operations amount to local transformations of (and in local dependence on) the tape contents.
And computational cost is measured in dependence on the binary length $n\in\IN$ of the input $\vec x\in\TWO^n$.
\item[b)]
Computing over the countable space $X$ of finite graphs proceeds by encoding them over $\TWO^*$:
for example as adjacency matrices' binary entries.
Cost parameter $n\in\IN$ commonly denotes the number of nodes of the graph, or the binary length of the encoded matrix:
both are polynomially related to each other.
\item[c)]
Consider the countable space $X=\IN$ of natural numbers.
Turing-computation encodes $\IN$ also over $\TWO^*$ from (a),
for instance in binary or in unary:
These are computably related, and therefore induce the same notions of computability on $\IN$;
but their lengths, and induced notions of complexity, are not polynomially related.
\item[d)]
Recall the type-2 machine model \cite[\S2.1]{Wei00} operating on 
the Cantor space $\Cantor:=\TWO^\IN$ of \emph{in}finite binary sequences
equipped with the metric
$D:(\bar u,\bar v)\;\mapsto\;2^{-\min\{n:u_n\neq v_n\}}$.
The computation of such a machine amounts to a partial transformation $F:\bar u\mapsto\bar v=(v_n)$ from/to $\Cantor$.
Equivalently via currying, the machine on input of `infinite' $\bar u\in\Cantor$ and of `finite' $\sdone^n$ (encoding $n\in\IN$ in unary)
prints $v_n$: the $n$-th entry of the output sequence $\bar v\in\Cantor$.
\\
Computational cost 
is commonly gauged in dependence of the index position $n$ within the binary \emph{out}put sequence,
that is, the length of the finite initial segment $\bar v_{<n}=(v_0,v_1,\ldots,v_{n-1})$ written so far \cite[\S7.1]{Wei00};
equivalently: the length $n$ of the finite part $\sdone^n$ of the \emph{in}put.
In either case, a computational cost bound $t=t(n)$ must not depend on the \emph{in}finite part $\bar u$ of the input.
It follows that any partial $F:\subseteq\Cantor\to\Cantor$ computed in time $t:\IN\to\IN$ has modulus of continuity $t$
according to Subsection~\ref{ss:Metric}.
\item[e)]
Computation on the real unit interval $X=[0;1]$ proceeds by encoding $X$ over ground space $\Cantor$, in various ways:
\item[e\,i)]
The binary expansion renders tripling $[0;1/3]\ni r\mapsto 3r\in[0;1]$ 
uncomputable
and is therefore unsuitable \cite[p.546]{Tur38}.
To see this, start with $1/6=(\sdzero.\sdzero \sdzero \sdone \sdzero \sdone\sdzero\sdone\ldots)_2$ and
consider, arbitrarily far off in the expansion, 
either replacing some single digit $\sdzero$ with $\sdone$
or replacing some single digit $\sdone$ with $\sdzero$:
The result of tripling in the first case will begin with
$(\sdzero.\sdone\ldots)_2$ 
yet in the second case begin with 
$(\sdzero.\sdzero\ldots)_2$.
To algorithmically output the first binary digit after the binary point
thus would require `knowing' all infinitely many binary digits of the input:
contradiction \cite[Exercise~7.2.7]{Wei00}.
\item[e\,ii)]
Let $r\in[0;1]$ be encoded as any sequence
\begin{equation} \label{e:Rho}
\big(\bin(a_0),\bin(c_0),\ldots,\bin(a_n),\bin(c_n),\ldots\big) \;\in\;\Cantor
\end{equation}
of binary numerators and denominators $a_n,c_n\in\IN$ 
of rational approximations satisfying $|r-a_n/c_n|\leq1/2^n$. 
This representation induces the `standard' notion of real computability \cite{Grz57}. 
In particular every (qualitatively) computable partial real function must be (qualitatively) continuous \cite[Theorem~4.3.1]{Wei00}. 
But this representation induces no reasonable notion of computational complexity. \\
To see this, consider the number $t\in\IN$ of steps some algorithm makes;
and these steps include skipping some number $\mu(n)\leq t(n)$ of bits in the
infinite string from Equation~\eqref{e:Rho} 
in order to reach and extract $a_n/c_n$ as guaranteed approximation to $r$ up to $1/2^n$.
For $r:=1/2$ with $a_n:=1+2^{k_n}$ and $c_n:=2^{1+k_n}$ with sufficiently large integers $k_n\geq n$,
this number $\mu(n)$ can be arbitrarily large, meaning that the running time $t\geq\mu$ of
such an algorithm cannot be bounded in terms of the precision parameter $n$ \cite[Examples~7.2.1+7.2.3]{Wei00}. 
Note that, compared to the required approximation error bound $1/2^n$,  denominators $2^{k_n}$ in this counterexample are exceedingly large---which Item~(e\,iii) prevents.
\end{enumerate}\noindent
Here and in the sequel, $\bin:\IN\to\TWO^*$ denotes binary expansion of natural numbers with delimiter quelling leading $\sdzero$s:
\begin{multline*}
\bin:\IN\to\{\sdzero,\sdone\}^*, \quad 0\;\mapsto\;\sdzero\,\sdone, \\
\{2^n,\ldots,2^{n+1}-1\}\;\ni\; 2^n+2^{n-1}b_{n-1}+2^{n-2}b_{n-2}+\cdots+2b_1+b_0 \;\mapsto\; \\
\mapsto\; 
(\sdone\,\sdzero\,b_{n-1}\,\sdzero\,\ldots b_1\,\sdzero\,b_0\,\sdone)\; 
\in\; \{\sdzero,\sdone\}^{2n+2}  \enspace .
\end{multline*}
\begin{enumerate}
\item[e\,iii)]
Consider encoding $r\in[0;1]$ over $\Cantor$ as sequence $\big(\bin(a_1),\ldots,\bin(a_n),\ldots\big)$ 
of numerators $a_n\in\IN$ of dyadic approximations satisfying $|r-a_n/2^n|\leq1/2^n$. 
This induces the same notion of real computability as (e\,ii).
And, again, computational cost counts the number $t(n)$ of steps until output of (the numerator $a_n$ to) an approximation up to absolute error $1/2^n$;
or until output of the $n$-th bit of $\big(\bin(a_1),\ldots,\bin(a_m),\ldots\big)\in\Cantor$:
Both parameters are polynomially related to each other: 
$\big(\bin(a_1),\ldots,\bin(a_n)\big)$ has binary length between $n$ and $\calO(n^2)$, and in particular is bounded polynomially in $n$.
This renders arithmetic (addition, subtraction, multiplication, division) of real numbers computable in polynomial time.
Moreover any partial real function computed in time $t:\IN\to\IN$ has modulus of continuity $\poly\Big(t\big(\poly(n)\big)\Big)$;
cmp. Corollary~\ref{c:Main1} below.
\item[e\,iv)] 
Encode $x\in[0;1]$ as \emph{signed} binary expansion
$x=1/2+\sum_{m\geq0} s_m\cdot2^{-m-2}$ with $s_n\in\THREE\cong\{\sdzero\sdzero,\sdzero\sdone,\sdone\sdzero\}$.
Again, computational cost counts the number of steps until output of approximation $1/2+\sum_{m\leq n} s_m\cdot2^{-m-2}$
up to absolute error $1/2^n$; or until output of the $n$-th (signed) output bit of
$(s_0,\ldots,s_n,\ldots\big)$: Both are linearly related to each other.
They render addition of real numbers computable by a finite transducer
(see Figure~\ref{f:Transducer} below) and in particular 
in linear time \cite[Theorem~7.3.1]{Wei00}.
Note that $[0;1]^2\ni(x,y)\mapsto x+y\in[0;2]$ has modulus of continuity $\mu(n)=n+\calO(1)$.
Conversely, any partial real function computed in time $t$ has a modulus of continuity $\calO\Big(t\big(\calO(n)\big)\Big)$;
cmp. Corollary~\ref{c:Main1} below.
\item[f)]
Recall the oracle Turing machine model, 
allowed to query some oracle\footnote{We here deliberately restrict to `classical', Boolean-valued oracles.} 
$\varphi\subseteq\TWO^*$.
Equivalently via currying: query a total string predicate $\varphi:\TWO^*\to\TWO$.
Considering the oracle as variable \cite[\S3]{KC12} leads to computation
as a transformation $\Lambda:\varphi\mapsto\psi$ from/to the space $\Baire=\dual{(\TWO^*)}$ 
of total string predicates;
equivalently via currying: $\Lambda:\dual{(\TWO^*)}\times\TWO^*\ni(\varphi,\vec x)\mapsto \psi(\vec x)\in\TWO$.
\\
Cost here is to be bounded in dependence of the binary input length $n=|\vec x|$ of the finite part of the input,
but independently of the infinite part $\varphi\in\dom(\Lambda)$.
Equipping $\dual{(\TWO^*)}$ with the metric
$D:(\varphi,\psi)\mapsto2^{-\min\{|\vec x|:\varphi(\vec x)\neq\psi(\vec x)\}}$,
any partial $\Lambda:\subseteq\dual{(\TWO^*)}\to\dual{(\TWO^*)}$ computable in time $t$
has modulus of continuity $t$.
\item[g)]
Spaces $X$ of continuum cardinality beyond real numbers (e) are also commonly encoded
over Cantor space $\Cantor=\TWO^\IN$ \cite[\S3]{Wei00}, or over oracle space $\Baire=\dual{(\TWO^*)}$ \cite[\S3.4]{KC12};
see \cite{ZZ99,Zie03,Sch07} for some examples.
\\
Computability-theoretically `reasonable' representations are the \emph{admissible} ones
\cite[Theorem~3.2.9.2]{Wei00}; cmp. \cite{KW85,Sch02}:
(Precisely) these make continuous functions correspond to continuous translations on Cantor/oracle space and vice versa
\cite[Theorem~3.2.11]{Wei00}.
The encodings of real numbers from Items~(e\,ii) to (e\,iv) are admissible, but the binary one in (e\,i) is not.
\item[h)] Search problems admit possibly more than one answer/output to a given input.
For example the Fundamental Theorem of Algebra assigns to any $d$-tuple of complex coefficients $(c_0,\ldots,c_{d-1})\in\IC^d$
\emph{some} $d$-tuple $(z_1,\ldots,z_d)\in\IC^d$ of complex roots (including multiplicity)
of the polynomial $c_0+c_1Z+\cdots+c_{d-1}Z^{d-1}+Z^d$.
This non-uniqueness of the answer is modeled by a (set-valued aka) \emph{multi-}function
$\IC^d\toto\IC^d$, $(c_0,\ldots,c_{d-1})\mapstoto(z_1,\ldots,z_d)$; cmp. \cite{Luc77}.
\end{enumerate}
\end{example}
To summarize, common models of computation (Turing machine, type-2 machine, oracle machine)
naturally operate on certain `basic' spaces (finite binary strings $\TWO^*$, Cantor space $\Cantor=\TWO^\IN$, oracle space $\dual{(\TWO^*)}$)
with computational cost bounded in dependence of the binary length $n$ of the finite part of inputs,
but independently of the infinite parts.
Computation on other spaces proceeds by encoding them over such basic spaces;
and computing a function between such spaces amounts to translating codes of arguments to codes of values.

Note that, for computability questions,
\emph{equilogical} spaces have been suggested as basic spaces \cite{BBS04}.
Regarding complexity questions, we suggest to generalize $\Cantor$ and $\dual{(\TWO^*)}$,
namely to consider compact ultrametric spaces of diameter 1 as basic spaces
over which to encode other spaces. This choice will be justified in Subsection~\ref{ss:Selection}.
Recall that an ultrametric $D$ satisfies the following \emph{strong} triangle inequality:
\begin{equation}
\label{e:Ultrametric}
D(x,z) \;\leq\; \max\big\{\:D(x,y)\:,\:D(y,z)\;\} \enspace .
\end{equation}
The present work addresses the following question:

\begin{question}
\label{q:Main}
Fix some compact metric space $X$ and a compact ultrametric ground space $\calU$.
Which encodings of $X$ over $\calU$ are `suitable' with respect to computational complexity?
\end{question}
Recall (Example~\ref{x:Example}g) that, for computability (rather than complexity) purposes,
\emph{admissibility} has been established as answer to this Question~\ref{q:Main}.

\subsection{Previous/Related Work and Overview}
Computability, and qualitative coding theory, for continuous data 
dates back at least to Turing for the real case \cite{Tur37};
and to Kreitz and Weihrauch for the case of second-countable T0 spaces \cite{KW85},
see Fact~\ref{f:KreitzWeihrauch} below and the standard textbook \cite{Wei00};
with generalization to further spaces mostly due to Schr\"{o}der \cite{Sch01,Sch02,DBLP:conf/cie/Schroder06}.

Early rigorous bit-complexity investigations, and quantitative coding considerations over the reals
\cite{KF82,Fri84,Mue86a,DK89}, had culminated in the textbook \cite{Ko91};
see \cite{Kaw10,KC12,KO14} for recent work---and 
\cite{Wei03,KMRZ15,DBLP:conf/lics/SchroderS17,DBLP:journals/lmcs/Steinberg17}
for first algorithmic cost and coding considerations of spaces beyond continuous real functions.

The bit-complexity questions pursued there and in the present work 
are not to be confused with Information-Based/Algebraic Complexity Theory based on the 
unit-cost aka Blum-Shub-Smale aka real-RAM model \cite{IBC,BCSS,ACT,compGeom}.
The latter model also underlies the Nyquist-Shannon Sampling Theorem,
revolving around encoding certain real functions by the least number of real numbers---each one still an infinite object:
The present work is closer to Approximation and Information Theory \cite{ApproximationTheory}.

Embedding \emph{finite/countable} metric spaces into trees
(as particular cases of ultrametric spaces) with least distortion
is a well-studied problem \cite{DBLP:reference/cg/IndykM04};
and the present work considers encoding uncountably infinite spaces 
over the infinite tree of Cantor space.

\medskip
Subsection~\ref{ss:Metric} collects some mathematical basics of metric spaces,
with emphasis on quantitative properties like entropy and modulus of continuity.
Subsection~\ref{ss:Representations} recalls the Kreitz-Weihrauch framework
of encoding continuous data over Cantor space. Their qualitative key property,
\emph{admissibility}, yields the \emph{Main Theorem} (Fact~\ref{f:KreitzWeihrauch}).
Four ways of encoding the real unit interval illustrate important differences regarding computability and complexity.
They serve as guide in Section~\ref{s:Admissible2} to our new quantitative version of admissibility (Definition~\ref{d:Admissible2});
and yield the new, quantitative \emph{Main Theorem}~\ref{t:Main1} in Subsection~\ref{ss:Main1},
proven in Subsection~\ref{ss:Proofs1}. Theorem~\ref{t:Standard} asserts a large class
of compact metric spaces to admit quantitatively admissible representations over Cantor space,
including a generalization of the signed binary representation.

\begin{table}[htb]\normalsize
\begin{tabular}{c||cccccc}
\textbf{qualitative}  & computability & topology & (uniform) continuity & compactness  & equilogical \\[0.5ex] \hline
\textbf{quantitative} &  complexity &  metric     &  modulus of continuity &  entropy     &  ultrametric
\end{tabular}
\caption{\label{t:Dictionary}``Dictionary'' between
qualitative and quantitative coding of continuous data.}
\end{table}

Section~\ref{s:Multifunc} takes an alternative approach to quantitatively refine qualitative admissibility:
from the perspective of a logical (as opposed to topological) 
notion of continuity for \emph{multi}functions from \cite{BrHe94,PZ13}.
To this end Subsection~\ref{ss:Multifunc} recalls some terminology concerning multifunctions.
Then Subsection~\ref{ss:Continuity} motivates and introduces a quantitative version of continuity for multifunctions.
Our quantitative selection theorem for multifunctions between compact ultrametric spaces
in Subsection~\ref{ss:Selection} applies and builds on this notion. 
It leads in Subsection~\ref{ss:Ultra} to a quantitative \emph{Main Theorem}~\ref{t:Main2} 
for multifunctions between spaces equipped with quantitatively admissible representations
over suitable generalized (namely quantitatively homogeneous and quantitatively admissible) domains.

\subsection{Recap of Metric Spaces}
\label{ss:Metric}

Like resource-bounded complexity quantitatively refines qualitative computability, this subsection recalls relevant
metric properties that quantitatively refine qualitative topological properties of the mathematical spaces under consideration.

Following {\rm\cite[\S6]{Wei03}}, a (binary) modulus 
of continuity of a function $f:X\to Y$ between metric spaces $(X,d)$ and $(Y,e)$ 
is a non-decreasing mapping $\mu:\IN\to\IN$ such that
\begin{equation}
\label{e:Modulus}
d(x,x')\;\leq\;2^{-\mu(n)} \quad\Rightarrow\quad e\big(f(x),f(x')\big)\;\leq\;2^{-n} \enspace .
\end{equation}
Equivalently: the \emph{induced semi-metric}
\begin{equation}
\label{e:Modulus2}
e_f:f[X]\times f[X]\;\ni\;(y,y')\;\mapsto\; \min\big\{ d(x,x') \::\: f(x)=y, f(x')=y' \big\}
\end{equation}
satisfies $e_f(y,y')\leq2^{-\mu(n)}\Rightarrow e(y,y')\leq2^{-n}$.
Note that $e_f$ is well-defined whenever $f$ is continuous with compact domain,
but may violate the triangle inequality. 
Call a modulus of continuity $\mu$ of $f$ \emph{minimal} if it is pointwise minimal 
among those satisfying Equation~\eqref{e:Modulus} and/or \eqref{e:Modulus2}.

\begin{figure}[htb]
\begin{center}\includegraphics[width=0.8\textwidth]{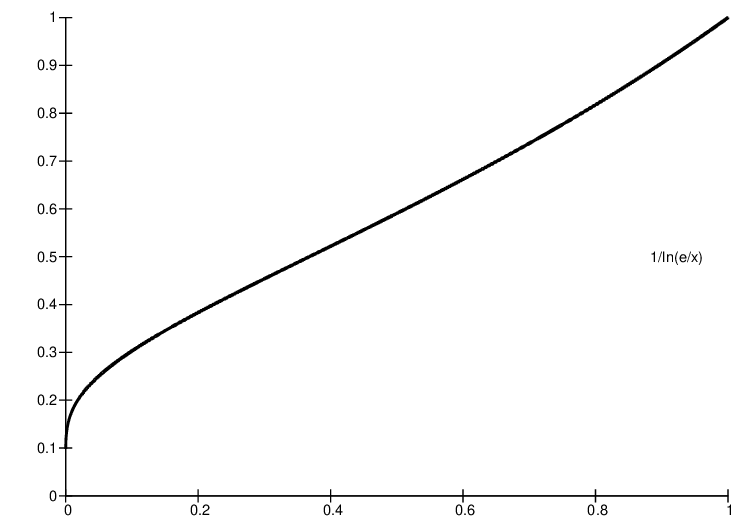}\end{center}
\caption{\label{f:superpol}}
\end{figure}

\begin{remark}
\label{r:Modulus}
The integer $m=\mu(n)$ roughly captures the number of bits of approximations 
to arguments $x$ sufficient to guarantee $n$ bits of approximation to $y=f(x)$.
\begin{enumerate}
\item[a)]
Precisely the uniformly continuous functions have a (minimal) modulus of continuity.
\item[b)]
Suppose $X$ is compact or convex complete or $Y$ is bounded.
Then Lipschitz-continuous functions $f:X\to Y$ are precisely those with modulus $\mu(n)=n+\calO(1)$.
And H\"{o}lder-continuous functions $f:X\to Y$ are precisely those with linear modulus $\mu(n)=\calO(n)$.
\item[c)]
If $f:X\to Y$ has modulus of continuity $\mu$ and $g:Y\to Z$ has modulus of continuity $\nu$,
the $g\circ f$ has modulus of continuity $\mu\circ\nu$.
\item[d)]
The function $h:(0;1]\ni t\mapsto 1/\ln(e/t)\in(0;1]$
from Figure~\ref{f:superpol}
extends uniquely continuously to $0$.
It has an exponential, but no polynomial,
modulus of continuity. 
\item[e)]
The $k$-fold iterate of $h:[0;1]\to[0;1]$ from (d)
has modulus of continuity an exponential tower of height $k$,
but not of height $k-1$.
\end{enumerate}
\end{remark}
For self-containment, a proof is provided at the end of this subsection.
Recall that a (not necessarily linear) metric space $(X,d)$ is \emph{convex} if,
to any $x,z\in X$ there exists $y\in X$ distinct from $x,z$ with $d(x,z)=d(x,y)+d(y,z)$. 
Abbreviate with $\cball(x,r)=\{x':d(x,x')\leq r\}\subseteq X$ 
the closed ball with center $x\in X$ and radius $r>0$,
and let $\diam(X)=\sup\{d(x,y):x,y\in Y\}$ denote the diameter of $X$.
The present work considers representations over compact ultrametric spaces of diameter 1,
a rich class of spaces:
\begin{lemma}
\label{l:Ultrametric}
\begin{enumerate}
\item[a)]
Any finite set, when equipped with the discrete metric, constitutes a compact ultrametric space of diameter 1. \\
\item[b)]
Equip the space $\TWO^*$ of finite binary strings with the ultrametric
\[ D(\vec u,\vec v) \;:=\; 2^{-\min\{n\::\:n\geq|\vec u| \:\vee\: n \:\geq\: |\vec v| \:\vee\: u_n\neq v_n\}} \enspace : \]
Its completion $\TWO^*\uplus\Cantor$ is compact with diameter 1. 
Same for the space $\TWO^+$ of \emph{non-}empty finite binary strings with ultrametric $2D$.
\item[c)]
Let $(U,D)$ and $(V,E)$ denote compact (ultra)metric spaces of diameter 1.
The Cartesian product $U\times V$ is again an (ultra)metric space of diameter 1
when equipped with the maximum metric 
\begin{equation}
\label{e:MaxMetric}
D\times E \;:=\; \max\{D,E\} \;:\; \big((u,v),(u',v')\big) \;\mapsto\; \max\big\{D(u,u'),E(v,v')\big\} \enspace .
\end{equation}
Similarly/inductively for the Cartesian product of finitely many (ultra)metric spaces of diameter 1.
\item[d)]
Let $(U,D)$ denote a compact metric space and $(V,E)$ a compact (ultra)metric space of diameter 1.
Consider the space $\LIP(U,V)$ of total functions $f:U\to V$
which are non-expansive (aka 1-Lipschitz) in satisfying
\begin{equation}
\label{e:Nonexpansive}
\forall u,u'\in U: \quad E\big(f(u),f(u')\big) \;\leq\; D(u,u') \enspace .
\end{equation}
Equipped with the supremum metric $E^U(f,g)=\sup\big\{e\big(f(u),g(u)\big):u\in U\big\}$,
$\LIP(U,V)$ is again a compact (ultra)metric space of diameter 1.
\item[e)]  
Fix compact (ultra)metric spaces $(U_j,D_j)$ of diameter 1, $j\in J\subseteq\IN$. 
Their Cartesian product $\prod_{j\in J} U_j$ is again compact of diameter $\leq1$ when equipped with the (ultra)metric 
\begin{equation}
\label{e:HilbertMetric}
\prod\nolimits_{j\in J} D_j/2^j \::\: \big((u_j)_{_j},(u'_j)_{_j}\big) \;\mapsto\; \max\big\{ D_j(u_j,u'_j)/2^{j} : j\in J \big\} \enspace .
\end{equation} 
\end{enumerate}
\end{lemma}
Note that Cantor space $\Cantor$ is recovered from Lemma~\ref{l:Ultrametric}e) with $U_j:\equiv\TWO$ according to (a).
Also oracle space $\dual{(\TWO^*)}$ arises from Lemma~\ref{l:Ultrametric}e), 
namely with the finite Cartesian products $U_j:=\prod_{\TWO^j} \TWO=\TWO^{\TWO^j}$, $j\geq0$, equipped with the discrete metric according to (a). 
$\LIP(\Cantor,\Cantor)$ is a compact ultrametric space of diameter 1 according to Lemma~\ref{l:Ultrametric}d);
and isometric to $\dual(\TWO^+)=\prod_{j\geq0} \TWO^{(\TWO^{j+1})}$ according to Lemma~\ref{l:Ultrametric}e). 
The following lemma generalizes this isometry:

\begin{lemma}
\label{l:Donghyun}
Let $U_j,V_j$ denote sequences of finite sets, equipped with the discrete metric.
Consider the compact ultrametric spaces $U=\prod_j U_j$ and $V=\prod_j V_j$
according to Lemma~\ref{l:Ultrametric}e).
Then $\LIP(U,V)$ according to Lemma~\ref{l:Ultrametric}d) is isometric to
$\prod_j W_j$ according to Lemma~\ref{l:Ultrametric}e), where
\[ W_j \;:=\; \prod\nolimits_{U_0\times U_1\times\cdots\times U_j} V_j \;=\; V_j^{U_0\times U_1\times\cdots\times U_j} \enspace . \]
\end{lemma}
Recall \cite{Kolmogorov} that the \emph{unary entropy} of $X$,
also known as its \emph{width} \cite[\S6]{Wei03},
is the mapping $\unary{\Entropy}=\unary{\Entropy}_X:\IN\to\IN\cup\{\infty\}$ such that
$X$ can by covered by $\unary{\Entropy}(n)$, but not by less,
closed balls $\cball(x,2^{-n})$ of radius $2^{-n}$ around centers $x\in X$.
To match the binary conception of the modulus \eqref{e:Modulus},
we here focus on the \emph{binary} entropy 
$\Entropy_X=\lceil\log_2\unary{\Entropy}_X\rceil$,
and call it just entropy \cite{Walters,DBLP:conf/lics/KawamuraS016,DBLP:journals/lmcs/Steinberg17}; 
cmp. Remark~\ref{r:Skolem} below.

\begin{figure}[htb]
\begin{center}\includegraphics[width=0.95\textwidth]{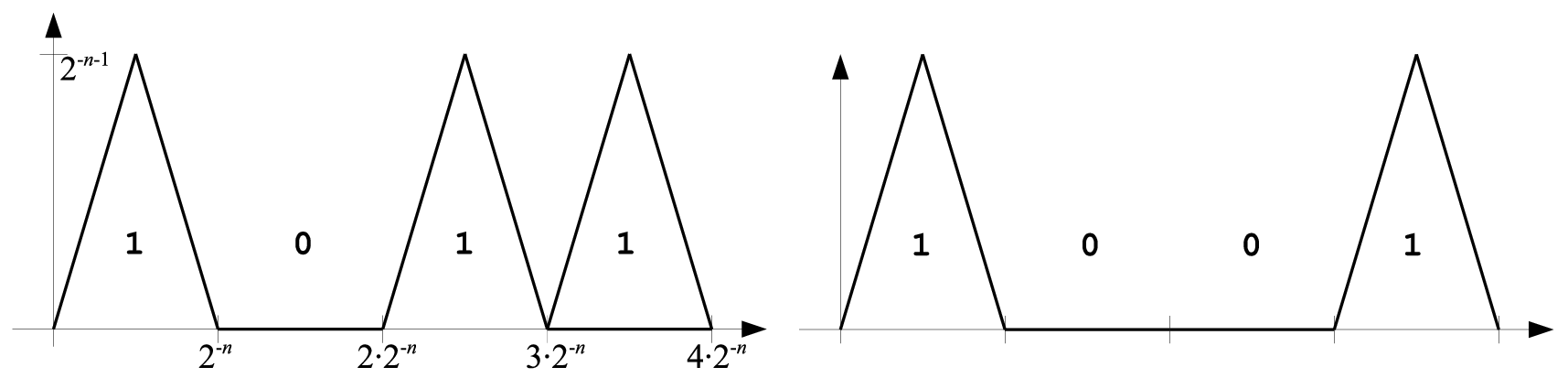}\end{center}
\caption{\label{f:manyfunc}}
\end{figure}

\begin{lemma}
\label{l:Entropy}
The integer $m=\Entropy(n)$ roughly describes the number of bits sufficient and necessary
to specify, for every possible $x\in X$, some approximation up to error $2^{-n}$.
\begin{enumerate}
\item[a)]
Cantor space has entropy $\Entropy(n)=n$;
the real unit interval $[0;1]$ has entropy $\Entropy(n)=\max\{n-1,0\}$.
The Hilbert Cube $\Hilbert=\prod_{j\geq0}[0;1]$ with metric (\eqref{e:HilbertMetric})
has quadratic entropy $\Entropy(n)=\Theta(n^2)$. 
\item[b)]
Generalizing Example~\ref{x:Example}f),  
fix some non-empty $\Phi\subseteq\TWO^*$ and let $\dual{\Phi}:=\TWO^{\Phi}$ denote
the collection of total string predicates $\varphi:\Phi\to\TWO$.
Equip $\dual{\Phi}$ with the ultrametric $D(\varphi,\varphi')=1/2^{\min{|\vec v|:\varphi(\vec v)\neq\varphi(\vec v')}}$.
$\dual{\Phi}$ is compact with entropy $\Entropy_{\dual{\Phi}}(n)=\Card\Phi_{<n}$, where $\Phi_{<m}:=\Phi\cap\TWO^{<m}$
and $\TWO^{<m}=\bigcup_{n<m}\TWO^n$.
\item[c)]
Let compact $(X,d)$ and $(Y,e)$ have entropies $\Entropy_X$ and $\Entropy_Y$, respectively.
Then the entropy $\Entropy_{X\times Y}$ of compact $\big(X\times Y,\max\{d,e\}\big)$ satisfies
\[ \forall n>0: \quad \Entropy_X(n-1)+\Entropy_Y(n-1) \;\leq\; \Entropy_{X\times Y}(n)+1 \;\leq\; \Entropy_X(n)+\Entropy_Y(n)+1 \enspace . \]
\item[d)]
Let compact $(X_j,d_j)$ have diameters $\leq1$ and entropy $\Entropy_j$ for each $j\in\IN$.
Then $\big(\prod_jX_j,\sup_j d_j/2^j\big)$ is compact and has entropy $\Entropy$ satisfying 
\[ \forall n>0: \quad \sum\nolimits_{j<n} \Entropy_j(n-j-1)  \;\leq\; 
1\:+\:\Entropy(n) \;\leq\; 1\:+\:\sum\nolimits_{j\leq n} \Entropy_j(n-j)  \enspace .  \]
\item[e)]
Fix a compact metric space $(X,d)$ with entropy $\Entropy$.
Let $\calK(X)$ denote the set of non-empty closed subsets of $X$
and equip it with the Hausdorff metric
\begin{equation}
\label{e:Hausdorff}
d_\calK(V,W) \;\mapsto\; \max\big\{\sup\{ d_V(w) :w\in W\},\sup\{d_W(v):v\in V\}\big\}, \end{equation}
where $d_V:X\ni x\mapsto\inf\{ d(x,v) : v\in V\}\in\LIP(X,\IR)$
denotes $V$'s distance function.
Then $\big(\calK(X),d_\calK\big)$ constitutes a compact metric space \cite[Exercise~8.1.10]{Wei00}.
It has entropy $H\leq2^{\Entropy}$ with $2^{\Entropy(n)-1}<H(n+1)$.
\item[f)]
Fix a connected compact metric space $(X,d)$ 
with entropy $\Entropy$,
and consider the convex metric space $\dual{X}:=\LIP(X,[0;1])$ of non-expansive (=1-Lipschitz) real functions 
equipped with the supremum norm $|f|=\sup_{x\in X} |f(x)|$, compact by Arzel\'a-Ascoli.
This has entropy $\dual{\Entropy}(n):=\Entropy_{\dual{X}}(n)=\Theta\big(2^{\Entropy(n\pm\calO(1))}\big)$; 
more precisely it holds: $2^{\Entropy(n-1)-1} \;<\;
\dual{\Entropy}(n) \;\leq\;\calO\big(2^{\Entropy(n+2)}\big)$;
see Figure~\ref{f:manyfunc}.
\item[g)] 
For $f:X\to Y$ a mapping between metric spaces with modulus of continuity $\mu$,
the image $f[X]\subseteq Y$ has entropy
$\Entropy_{f[X]}\leq \Entropy_X\circ\mu$; cmp. \cite[Lemma~3.1.13]{SteinbergDisse}.
\item[h)]
Every connected compact metric space $X$ has entropy at least linear $\Entropy(n)\geq n+\Omega(1)$.
\item[j)]
Fix an arbitrary non-decreasing unbounded mapping $\nu:\IN\to\IN$
and re-consider Cantor space $\Cantor$ but now equipped with 
$D_\nu:(\bar a,\bar b)\mapsto 2^{-\nu(\min\{m:a_m\neq b_m\})}\in[0;1]$.
This constitutes a metric, topologically equivalent to $D$,
but now inducing entropy $\Entropy_\nu=\loinv{\nu}:n\mapsto\min\{ m\::\: \nu(m)\geq n\}$.
\end{enumerate}
\end{lemma}
The notation $X\mapsto\dual{X}$ shared between Items~(b) and (f) is in no danger of confusion,
regarding that the former maps discrete to discrete spaces 
while the latter maps connected to connected spaces.
Item~c) follows from Lemma~\ref{l:Entropy2}a+b) below, 
Item~d) from Lemma~\ref{l:Entropy2}a+c). 
For Item~h) refer to \cite[Example~48]{Donghyun}.
Item~g) constitutes a quantitative refinement of the qualitative fact
that images of compact sets under continuous functions are again compact.

\begin{remark}
\label{r:Skolem}
Indeed, Mathematical Logic suggests \emph{Skolemization}
as generic quantitative refinement of qualitative $\Pi_2$ (\ie $\forall\exists$) properties
\cite[Ex~4.8.2+Def~17.106+p.285+\S15.4+p.379+Def~17.116]{Koh08a}.
Thus Skolemizing uniform continuity of $f:X\to Y$ 
\[ 
 \forall n\in\IN \; \exists m\in\IN \; \forall x,x'\in X: \quad d(x,x')\leq1/2^m \;\Rightarrow\; e(y,y')\leq1/2^n, \qquad
y:=f(x), \; y':=f(x') 
\] 
yields the above notion of \emph{modulus}. 
We consider modulus as a non-decreasing number-theoretic mapping with arguments and values as exponents to base $1/2$:
due to its close connection with asymptotic computational cost, recall Example~\ref{x:Example}d).
Similarly, the unary and binary \emph{entropies} in Lemmas~\ref{l:Entropy}+\ref{l:Entropy2} 
arise as Skolemizations of qualitative pre-compactness. 
\end{remark}

\begin{proof}[Proof of Remark~\ref{r:Modulus}]
In case $(X,d)$ is complete then $(X,d)$ is convex iff,
to any $x,z\in X$, there exists a (not necessarily unique) isometry $\gamma=\gamma_{x,z}:[0;d(x,z)]\subseteq\IR\to X$ with $\gamma(0)=x$ and $\gamma\big(d(x,z)\big)=z$.
\begin{enumerate}
\item[b)]
Suppose $f$ is $2^\ell$-Lipschitz: $e(y,y')\leq 2^\ell\cdot d(x,x')$, where $y:=f(x)$ and $y':=f(x')$. 
Then $f$ has modulus $\mu(n)=n+\ell$, since $d(x,x')\leq2^{-(n+\ell)}$ implies 
$e(y,y') \leq 2^{\ell}\cdot d(x,x') \leq 2^{\ell}\cdot2^{-(n+\ell)} =2^{-n}$.

For the converse, let $f$ have modulus $\mu(n)=n+\ell$ and first suppose $\diam(Y)\leq K$ for $K\geq 2$.
Then $2^{-(n+\ell+1)}<d(x,x')\leq2^{-(n+\ell)}$ implies 
\[ e(y,y')\;\leq\; 2^{-n} \;=\; 2^{\ell+1}\cdot2^{-(n+\ell+1)} \;<\; 2^{\ell+1}\cdot d(x,x') \;\leq\; K\cdot 2^{\ell}\cdot d(x,x') \]
for $n\in\IN$; while in case $d(x,x')>2^{-\ell}$ it holds
\[ e(y,y') \;\leq\; K \;=\; K\cdot 2^{\ell}\cdot 2^{-\ell} \;<\; K\cdot 2^{\ell}\cdot d(x,x') \enspace :\]
showing that $f$ is $K\cdot2^{\ell}$-Lipschitz.

If $X$ is compact, then so is its image $f[X]\subseteq Y$ under continuous $f$
and in particular bounded; w.l.o.g. $Y=f[X]$, and we are back in the first case.

Finally let $f$ have modulus $\mu(n)=n+\ell$ and suppose $X$ is convex and complete.
Then $2^{-(n+\ell+1)}<d(x,x')\leq2^{-(n+\ell)}$ implies 
\[ e(y,y')\;\leq\; 2^{-n} \;=\; 2^{\ell+1}\cdot2^{-(n+\ell+1)} \;<\; 2^{\ell+1}\cdot d(x,x') \;\leq\; (K+1)\cdot 2^{\ell}\cdot d(x,x') \]
for $n\in\IN$ and $K\geq1$ defined in a minute. 
On the other hand in case $r:=d(x,x')>2^{-\ell}$, as explained above,
there exists an isometry $\gamma:[0;r]\to X$ with $\gamma(0)=x$ and $\gamma(r)=x'$.
In particular $x_k:=\gamma(k/2^{\ell})$ for $k=0,1,\ldots,K:=\lfloor r/2^\ell\rfloor\geq 1$ satisfy
\[ d(x_k,x_{k+1})\leq2^{-\ell} \quad\text{ and }\quad 
d(x,x') = \sum\nolimits_{k=0}^K d(x_k,x_{k+1}) \] 
with the abbreviation $x_{K+1}:=x'$. Hence, by triangle inequality, 
\begin{eqnarray*} e(y,y') &\leq& \sum\nolimits_{k=0}^K e\big(f(x_k),f(x_{k+1})\big)  
\;\leq\; \sum\nolimits_{k=0}^K 2^{-0} \\
&=& K+1 \;=\; (K+1)\cdot 2^{\ell} / 2^{\ell} \;\leq\; (K+1)\cdot 2^{\ell} \cdot d(x,x') 
\end{eqnarray*}
showing that $f$ is $(K+1)\cdot2^{\ell}$-Lipschitz.

Suppose $f$ is H\"{o}lder continuous,
$e(y,y')\leq 2^\ell\cdot d(x,x')^{1/k}$ for $k,\ell\in\IN$.
Then $f$ has modulus $\mu(n)=k\cdot (n+\ell)$.
Conversely let $f$ have modulus $\mu(n)=k(n+\ell)$ and suppose $\diam(Y)\leq K$.
Then $2^{-k(n+\ell+1)}<d(x,x')\leq2^{-k(n+\ell)}$ implies 
\[ e(y,y')\;\leq\; 2^{-n} \;=\; 2^{\ell+1}\cdot(2^{-k(n+\ell+1)})^{1/k} \;<\; 2^{\ell+1}\cdot d(x,x')^{1/k}  \]
for $n\in\IN$; while in case $d(x,x')>2^{-k\ell}$ it holds
\[ e(y,y') \;\leq\; K \;=\; K\cdot 2^{\ell}\cdot 2^{-\ell} \;<\; K\cdot 2^{\ell}\cdot d(x,x')^{1/k} . \]
The remaining cases proceeds similarly.
\end{enumerate}\end{proof}

\subsection{Recap of Representations}
\label{ss:Representations}

Coding one space $X$ over another space $U$ is formalized as a \emph{representation}:
a surjective partial mapping $\xi:\subseteq U\twoheadrightarrow X$;
cmp. \cite[\S3]{Wei00} or \cite[\S3.4]{KC12}.
Every $u\in\dom(\xi)$ is called a \emph{name} of $x=\xi(u)\in X$;
and for another representation $\upsilon$ of $Y$,
a $(\xi,\upsilon)$-\emph{realizer} is a mapping $F:\dom(\xi)\to\dom(\upsilon)$
such that it holds $f\circ\xi\sqsubseteq\upsilon\circ F$:
$F$ transforms any $\xi$-name of $x\in\dom(f)$ to some $\upsilon$-name of $y=f(x)$.
Here $g\sqsubseteq h$ means that $g$ is a (not necessarily proper) restriction of $h$.
Examples~\ref{x:Binary},\ref{x:Rational},\ref{x:Dyadic},\ref{x:SignDigit} below
illustrate the effect of different representations of the real unit interval.
\emph{Admissibility} is a crucial criterion for a representation 
to be suitable for computability purposes \cite{KW85,Sch02}: 

\begin{definition}
\label{d:Admissible}
Fix topological spaces $X$ and $U$.
\begin{enumerate}
\item[a)]
For partial mappings $\xi,\xi':\subseteq U\to X$,
a \emph{continuous reduction from $\xi'$ to $\xi$} 
is a continuous mapping $G:\dom(\xi')\to\dom(\xi)$
with $\xi'\sqsubseteq\xi\circ G$. \\
In this case we write ``$\xi'\reduceq\xi$''.
\item[b)]
Call representation $\xi:\subseteq U\twoheadrightarrow X$ \emph{admissible} ~iff~ 
\item[i)] it is continuous, and 
\item[ii)] Every continuous surjective $\xi':\subseteq U\twoheadrightarrow X$ satisfies $\xi'\reduceq\xi$.
\item[c)]
Call $\xi$ \emph{standard} ~iff~ it satisfies (i') and (ii'):
\item[i')]
The topology on $X$ is the final one w.r.t. $\xi$: 
precisely the preimages $\xi^{-1}[Y]$ of open sets $Y\subseteq X$ are open in $\dom(\xi)$.
\item[ii')]
Every (\emph{not} necessarily surjective) continuous $\xi':\subseteq U\to X$ satisfies $\xi'\reduceq\xi$.
\end{enumerate}
\end{definition}
Note that every standard representation is admissible.
``$\reduceq$'' is a partial order; and admissible representations are those maximal among the continuous ones.
Standard (and hence admissible) representations exist for every T$_0$ space \cite[Theorem~3.2.9.2]{Wei00}.
Admissible representations are Cartesian closed; and yield the Kreitz-Weihrauch (aka \emph{Main}) Theorem:

\begin{fact}
\label{f:KreitzWeihrauch}
Standard (and hence admissible) representations exist (at least) 
for every second-countable T$_0$ space \cite[Theorem~3.2.9.2]{Wei00}.
\\
Fix admissible $\xi:\subseteq\Cantor\twoheadrightarrow X$ 
and $\upsilon:\subseteq \Cantor\twoheadrightarrow Y$.
A function $f:X\to Y$ is continuous ~iff~ it admits a continuous
$(\xi,\upsilon)$-realizer $F:\dom(\xi)\to\dom(\upsilon)$
\cite[Theorem~3.2.11]{Wei00}.
\end{fact}
It has been pointed out 
that admissibility and the Main Theory are more naturally stated
in terms of \emph{sequential} (rather than topological) continuity \cite{Sch01,Sch02,DBLP:conf/cie/Schroder06,Schroeder2020}.
Subsection~\ref{ss:Continuity} below suggests and investigates a notion of quantitative sequential continuity for \emph{multi}(valued) functions;
note that being a $(\xi,\upsilon)$-realizer of $f$ means being a selection of the multifunction $\upsilon^{-1}\circ f\circ\xi$.
Examples~\ref{x:Binary},\ref{x:Rational},\ref{x:Dyadic},\ref{x:SignDigit} below 
formalize the four encodings of the real unit interval from Example~\ref{x:Example}e)
as representations and collect their qualitative and quantitative, metric and computational properties:

\begin{example} \label{x:Binary}
The \emph{binary} representation of $X=[0;1]$ is the total mapping
\[ 
\beta:\Cantor
\;\ni\; \bar b \;\mapsto\; 
\sum\nolimits_{m=0}^{\infty} b_m2^{-m-1} \;{\in}\; [0;1] \enspace .
\] 
Note that $D(\bar b,\bar b')=2^{-n}$ implies
\[ \big|\beta(\bar b)-\beta(\bar b')\big|\;\leq\; \sum\nolimits_{m=n}^{\infty} 1\cdot 2^{-m-1} \;=\; 2^{-n}  \]
hence $\beta$ has modulus of continuity $\mu:=\id:n\mapsto n$, \ie, is 1-Lipschitz aka non-expansive.
\\
However $\beta$ is not admissible {\rm\cite[Theorem~4.1.13.6]{Wei00}};
$\beta$ does not admit a continuous realizer for instance of the continuous mapping $[0;1/3]\ni x\mapsto 3x\in[0;1]$;
cmp. {\rm\cite[Example~2.1.4.7+Exercise~7.2.7]{Wei00}}.
\end{example}
Note that every real number has either one or two binary expansions, \ie, $\beta$-names.
Intuitively, such uniqueness prevents admissibility:
Every real number has uncountably many different names
for the admissible representations in upcoming Examples~\ref{x:Rational}, \ref{x:Dyadic}, and \ref{x:SignDigit}.

\begin{example} \label{x:Rational}
The \emph{rational} representation of $[0;1]$ is the partial mapping
$\rho:\subseteq\Cantor\twoheadrightarrow[0;1]$ with
\begin{multline*}
\big(
\bin(a_0)\: \bin(c_0)\;
\bin(a_1)\: \bin(c_1)\;
\ldots
\bin(a_n)\: \bin(c_n)\;
\ldots \big)
\;\mapsto\; \lim\nolimits_j a_j/c_j , \\[0.5ex]
\dom(\rho) \:=\: \big\{
\big(  
\ldots\:
\bin(a_n)\:\bin(c_n)\;\:\ldots\big) \::\:
\exists x\in[0;1] \: |a_n/c_n-x|\leq 2^{-n}\big\} \enspace .
\end{multline*}
Representation $\rho$ is admissible \cite[Lemmas~4.1.4+4.1.6+Theorem~4.3.2]{Wei00}
and in particular continuous, but not uniformly continuous
(its domain is not compact) and thus has no modulus of continuity \cite[Examples~7.2.1+7.2.3]{Wei00}.
\end{example}
Intuitively, in the $2^{-n}$-approximation $a_n/c_n$,
numerator and denominator may be unnecessarily long.
In the following representation on the other hand,
the denominator of the $2^{-n}$-approximation is fixed to $2^n$,
and the numerator thus also bounded to have binary length $\calO(n)$:
\begin{example} \label{x:Dyadic}
Here we consider the binary encoding of natural numbers of \emph{given} length without delimiters:
\begin{equation}
\label{e:Binary2}
\bin_n:\{0,\ldots,2^n-1\}\;\ni\; b_0+2b_1+\cdots+2^{n-1}b_{n-1} \;\mapsto\; (b_0,\ldots,b_{n-1}) \;\in\; \{\sdzero,\sdone\}^n 
\end{equation}
The \emph{dyadic representation} of $[0;1]$ is the partial mapping
\begin{multline*}
\delta:\subseteq\Cantor\;\ni\;\big(
\bin_1(a_1)\:
\bin_2(a_2)\:
\ldots
\bin_n(a_n)\:
\ldots \big)
\;\mapsto\; \lim\nolimits_j a_j/2^j , \\[0.5ex]
\dom(\delta) \:=\: \big\{
\big(  
\ldots\:
\bin_n(a_n)\:\ldots\big) \::\:
2^n > a_n\in\IN, \;
|a_n/2^n-a_m/2^m|\leq 2^{-n}+2^{-m}\big\}
\end{multline*}
For instance $r=1/2$ has $\delta$-names 
$\big(\bin_1(1):\bin_2(2)\:\bin_3(4)\ldots\big)$
and $\big(\bin_1(0)\:\bin_2(1)\:\bin_3(3)\ldots\big)$ 
and $\big(\bin_1(1)\:\bin_2(3)\:\bin_3(5)\ldots\big)$
and many more.
\begin{enumerate}
\item[i)] $\delta$ has minimal quadratic modulus of continuity 
\[ \mu(n) \;=\; \sum\nolimits_{m=1}^{n} |\bin(a_m)| \;=\; \sum\nolimits_{m=1}^{n} m
\;=\; n\cdot (n+1)/2 \;\leq\; \calO(n^2) \enspace : \]
polynomial, but not linear.
\item[ii)] $\delta$ is admissible.
More precisely it satisfies the following quantitative
strengthening of Condition~(iii) in Definition~\ref{d:Admissible}:
To every (not necessarily surjective) partial function $\xi':\subseteq\Cantor\to [0;1]$
with modulus of continuity $\mu'$
there exists a mapping $G:\dom(\xi')\to\dom(\delta)$ 
with modulus of continuity $\Lambda:\mu(n)\mapsto\mu'(n)$ 
such that $\xi'=\delta\circ G$ holds.
\item[iii)]
To every $n\in\IN$ and every $r,r'\in[0;1]$ with $|r-r'|\leq2^{-n-1}$,
there exist $\delta$-names $\bar u$ and $\bar u'$ 
of $r=\delta(\bar u)$ and $r'=\delta(\bar u')$
with $D(\bar u,\bar u')\leq2^{-\mu(n)}$.
\end{enumerate}
Properties (i)--(iii) assert that $\delta$ is a \emph{polynomially standard} representation
in the sense of Definition~\ref{d:Admissible2} below.
\end{example}
Modulus $\Lambda:\mu(n)\mapsto\mu'(n)$ in (ii) means $\Lambda=\mu'\circ\loinv{\mu}$ 
with the semi-inverse $\loinv{\mu}(n):=\min\{ m\::\: \mu(m)\geq n\}$.
Note that $\loinv{\mu}:\IN\to\IN$ is totally defined when $\mu$ is unbounded, and satisfies 
\begin{equation}
\label{e:LoInv}
\loinv{\mu}\circ\mu\;\leq\;\id\;\leq\mu\circ\loinv{\mu} \enspace .
\end{equation}
For $r\in\IR$, let $\lfloor r\rceil$ denote the integer closest to $r$;
ties broken towards 0, \ie $\lfloor k+1/2\rceil=k$ and $\lfloor -k-1/2\rceil=-k$ for $k\in\IN$.
Note that rounding to nearest (rather than always down or always up)
\begin{equation}
\label{e:BetterApprox}
a_n \;:=\; \lfloor r\cdot2^n\rceil \quad\text{satisfies}\quad  |r-a_n/2^n|\leq2^{-n\pmb{-1}} 
\end{equation}
yields twice as good an approximation as the required error bound $2^{-n}$: crucial for (the proof of) Property~(ii).

\begin{proof}[Example~\ref{x:Dyadic}]
Let us abbreviate $\vec a_{<m}=(a_j:j<m)$ for $m\in\IN$ and an in/finite sequence $\vec a=(a_{j})_{_j}$.
\begin{enumerate}
\item[ii)] 
Fix $n\in\IN$ and consider the finite set $\vec\Xi_{m}$ of $\vec u\in\TWO^m$ 
with $(\vec u\circ\Cantor)\cap\dom(\xi')\neq\emptyset$,
that is, the set of all extensions of finite $\vec u$ to infinite strings in $\dom(\xi')$.
We record
\begin{gather}
\bigcup\nolimits _{\vec u\in\vec\Xi_m} \xi'[\vec u\circ\Cantor] \;=\; [0;1] \label{e:4Admissible1} \\
\vec\Xi_m \;=\; \big\{\vec u|_{<m}:\vec u\in\vec\Xi_{m+1}\big\} \label{e:4Admissible2} \\  
\diam\big(\xi'[\vec u\circ\Cantor]\big)\leq2^{-n}, \label{e:4Admissible3} \quad \text{hence} \\
r(\vec u)\;:=\; \max\big\{\xi'[\vec u\circ\Cantor]\big\}/2 \;+\; \min\big(\xi'[\vec u\circ\Cantor]\big)/2 \label{e:4Admissible4}
\end{gather}
has $\xi'[\vec u\circ\Cantor]\subseteq\cball\big(r(\vec u),2^{-n\pmb{-1}}\big)$ for all $\vec u\in\vec\Xi_m$ with $m:=\mu'(n)$. 
Therefore 
\[ 
F_n\;:\;\vec\Xi_{\mu'(n)}\;\ni\;\vec u\;\mapsto \;\lfloor 2^n\cdot r(\vec u)\rceil \;
\in\;\{0,\ldots,2^n-1\} \]
satisfies
$\xi'[\vec u\circ\Cantor] \subseteq \cball\big(F_{n}(\vec u)/2^n,2^{-n}\big)$
due to Equation~\eqref{e:BetterApprox}.
Finally 
\[ G:\dom(\xi')\to\dom(\delta), \quad \bar u \;\mapsto\; \Big( \bin_n\big(F_n(\bar u_{<\mu'(n)})\big) \;:\: n=1,2,\ldots \Big)  \]
maps $\dom(\xi')$ to $\dom(\delta)$,
and any $\xi'$-name of some $r\in[0;1]$ to a $\delta$-name of the same $r$: $\xi'\sqsubseteq\delta\circ G$.
Moreover, $G$ has modulus of continuity $\Lambda:\mu(n)\mapsto\mu'(n)$.
\item[iii)]
To $r\in [0;1]$ consider the $\delta$-name $\bar u:=\big(\ldots \bin_m(a_m)\:\ldots \big)\in\Cantor$ of $r$
with $a_m:=\lfloor r\cdot2^m\rceil$. 
In view of Equation~\eqref{e:BetterApprox}, it holds
$|r'-a_n/2^n|\leq|r'-r|+|r-a_n/2^n|\leq 2^{-n-1}+2^{-n-1}$; hence $(a_1,\ldots,a_n)$ can be extended 
to a $\delta$-name $\bar u'$ of $r'$ via $a'_{m}:=\lfloor r'\cdot2^{m}\rceil$ for all $m>n$.
Both agree on the initial segment $\big( \bin_1(a_1)\: \ldots \bin_m(a_n)\big)$ of binary length $\mu(n)$
and hence have $D(\bar u,\bar u')\leq2^{-\mu(n)}$. 
\qed\end{enumerate}\end{proof}
The dyadic representation has quadratic instead of linear 
modulus of continuity. This sub-optimality comes from the following redundancy: 
the $\calO(n)$ bits of $a_n$ providing approximation up to error $2^{-n}$
are preceded by the up to $\calO(n^2)$ bits of $a_1,\ldots,a_{n-1}$ providing worse approximations.
The \emph{signed} binary representation below on the other hand achieves both optimal modulus of continuity
(like the binary representation) and admissibility according to Definition~\ref{d:Admissible}.
It extends the standard binary digits $0$ and $1$ with $-1$:
\[ [0;1] \;\ni\; r \;=\; \tfrac{1}{2} \;+\; \sum\nolimits_m s_m2^{-m-2}, \quad s_m\in\THREE \] 
where $s_m\in\THREE$ is encoded as bit-tuples $\{\sdzero,\sdone\}^2\setminus \{\sdone\sdone\}$:
\begin{gather*}
\tilde\sigma:\THREE \;\ni\; s_m \;=\; 2b_{2m}+b_{2m+1}-1 \;\mapsto \; (b_{2m},b_{2m+1}) \;\in\; 
\{\sdzero\sdone, \sdone\sdzero, \sdzero\sdzero\}, \\
\qquad \THREE^\IN\ni(s_m)_{_m}\mapsto \big(\tilde\sigma(s_m)\big)_{_m}\in\TWO^\IN \enspace .
\end{gather*}

\begin{example}
\label{x:SignDigit}
The \emph{signed binary representation} of $[0;1]$ is the total mapping
\[ 
\sigma:\subseteq\;\Cantor
\;\ni\; \bar b \;\mapsto\; 
\tfrac{1}{2}\;+\;\sum\nolimits_{m=0}^{\infty} (2b_{2m}+b_{2m+1}-1) \cdot 2^{-m-2} 
\;{\in}\; [0;1] 
\] 
For instance $r=1/2$ has $\sigma$-names 
$\big(\bin(1):\bin(2)\:\bin(4)\ldots\big)$
and $\big(\bin(0)\:\bin(1)\:\bin(3)\ldots\big)$ 
and $\big(\bin(1)\:\bin(3)\:\bin(5)\ldots\big)$
and many more.
\begin{enumerate}
\item[i)]
$\sigma$ has linear modulus of continuity $\mu(n)=2n$.
\item[ii)] $\sigma$ is admissible \cite[Theorem~7.2.7+Subsection~7.3]{Wei00}. 
More precisely it satisfies the following quantitative
strengthening of Condition~(iii) in Definition~\ref{d:Admissible}:
To every (not necessarily surjective) partial function $\xi':\subseteq\Cantor\to [0;1]$
with modulus of continuity $\mu'$
there exists a mapping $G:\dom(\xi')\to\dom(\sigma)$ 
with modulus of continuity $\Lambda:2m\mapsto\mu'(m+1)$
such that $\xi'=\sigma\circ G$ holds.
\item[iii)]
To every $n\in\IN$ and every $r,r'\in[0;1]$ with $|r-r'|\leq2^{-n-1}$,
there exist $\sigma$-names $\bar b$ and $\bar b'$ 
of $r=\sigma(\bar b)$ and $r'=\sigma(\bar b')$
with $D(\bar b,\bar b')\leq2^{-2n}$.
\item[iv)] With respect to the signed binary representation, real averaging 
\[ [0;1]\times[0;1] \;\ni\; (x,y) \;\mapsto\; (x+y)/2 \;\in\;[0;1] \]
is computable in linear time, and in fact by a linear transducer.
\end{enumerate}
Properties (i)--(iii) assert that $\sigma$ is a \emph{linearly standard} representation
in the sense of Definition~\ref{d:Admissible2} below.
\end{example}
Intuitively speaking, appending more digits to some finite binary expansion of a real number
can only increase but not decrease its value; while 
Properties (ii) and (iii) of the signed binary expansion build on its ability
to improve any initial approximation in both directions, up and down.

\begin{figure}[htb]
\includegraphics[width=0.95\textwidth]{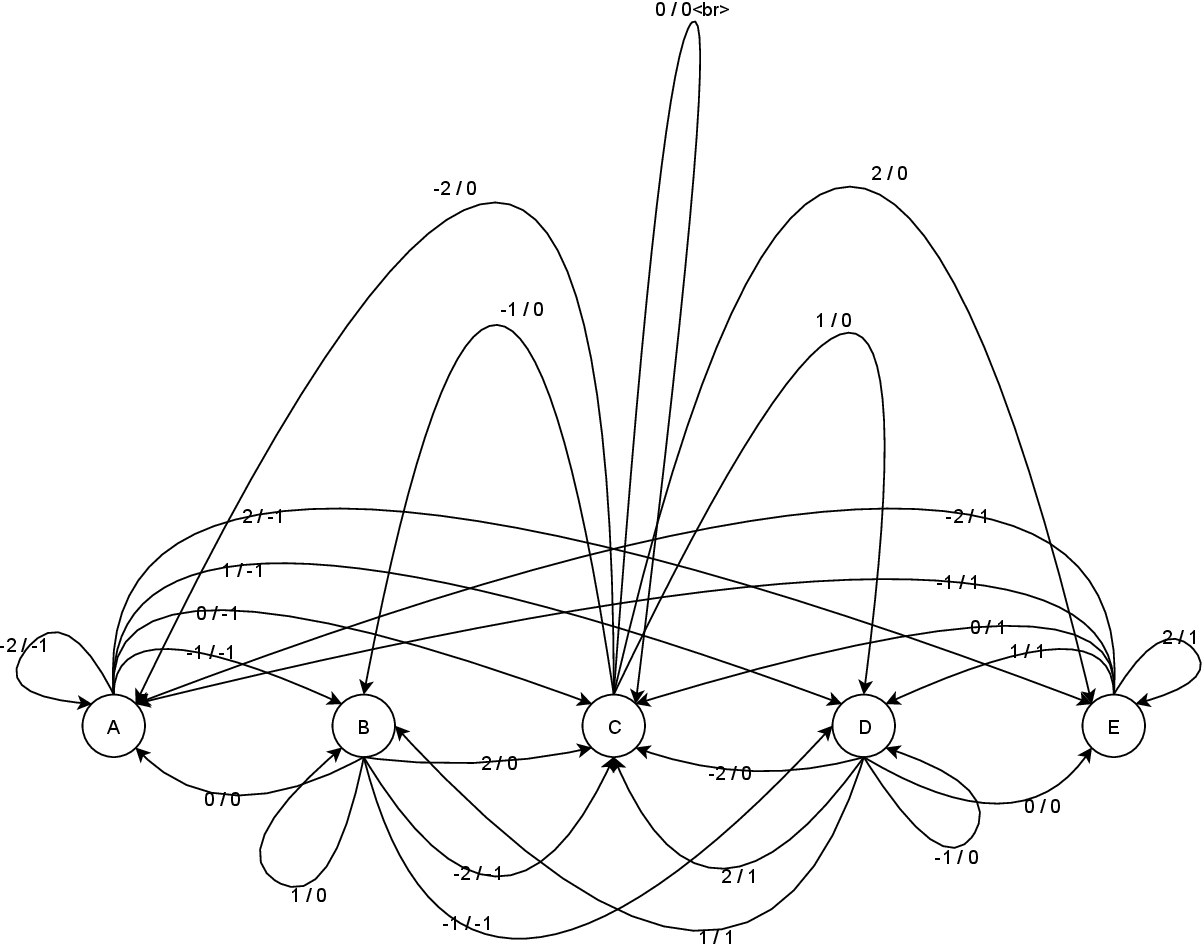}
\caption{\label{f:Transducer}Finite transducer
computing real averaging $(x,y)\mapsto(x+y)/2$ 
with respect to the signed-digit representation.
Edges are labeled with pairs: 
\emph{sum of both next signed-digit inputs}/\emph{signed-digit output}.}
\end{figure}

\begin{proof}
Let us abbreviate
\begin{gather*}
\tildesbin \;:\; \THREE^n \;\ni\; (s_0,\ldots,s_{n-1})
\;\mapsto\; \sum\nolimits_{m<n} s_m \cdot 2^{-m-1} \;\in\; [-1;1] \enspace , \\
\sbin:\THREE^* \;\ni\; \vec{s}\;\mapsto \; \tfrac{1}{2}\;+\; \tildesbin(\vec{s})/2 \;\in\; [0;1] \enspace , \\
\sbin:\THREE^\IN\;\ni\;\bar s\;\mapsto\;\lim_n \sbin\big(\bar s_{<n}\big)\;\in\;[0;1] \enspace ,
\end{gather*}
such that $\sigma=\sbin\circ\tilde\sigma$.
\begin{enumerate}
\item[i)] 
Regarding quantitative continuity,
changing $(b_{2m},b_{2m+1})=\sdzero\sdzero=\tilde\sigma(-1)$
to $(b'_{2m},b'_{2m+1})=\sdone\sdzero=\tilde\sigma(+1)$ for all $m\geq n$
changes $r=\sigma(\bar b)$ to $r'=\sigma(\bar b')$ with
$|r'-r|=\sum\nolimits_{m\geq n} 2 \cdot 2^{-m-2}=2^{-n}$
while $D(\bar b',\bar b)=2^{-2n}$.
Regarding the other case, consider 
changing $(b_{2n},b_{2n+1})=\sdzero\sdzero=\tilde\sigma(-1)$
to $(b'_{2n},b'_{2n+1})=\sdzero\sdone=\tilde\sigma(0)$ 
and $(b_{2m},b_{2m+1})=\sdzero\sdzero=\tilde\sigma(-1)$
to $(b'_{2m},b'_{2m+1})=\sdone\sdzero=\tilde\sigma(+1)$ for all $m>n$:
This changes $r$ to $r'$ with
$|r'-r|=2^{-n-2}+\sum\nolimits_{m>n} 2 \cdot 2^{-m-2}=3\cdot2^{-n-2}$
while $D(\bar b',\bar b)=2^{-2n-1}$.
Together it follows that $\mu(n)=2n$ is a modulus of continuity of $\sigma$.
\\
Regarding surjectivity of $\tildesbin$ (and in consequence of $\sbin$ and $\sigma$), we record (also for later use):
\item[i')] $\displaystyle 
\forall r\in[-1;1] \;
\forall n\in\IN \;  \exists \vec s\in\THREE^n: \quad |r \;-\; \tildesbin(\vec s)| \;\leq\; 2^{-n}$. \\
To see this in case $r\in[0;1)$, consider
the \emph{un}signed binary expansion of $\lfloor 2^{n}\cdot r\rfloor\in\big[0;2^{-n}big)$;
and negate it in case $r\in(-1;0]$. Note that least-significant signed digit $s_{n-1}$ has weight $2^{-n}$.
For $\sbin(\vec s)$ instead of $\tildesbin$, the error bound improves from $2^{-n}$ to $2^{-n-1}$.
\item[ii)] 
We re-use the notations $\vec\Xi_{m}$ and $r(\vec u)$ with
Properties~\eqref{e:4Admissible1},\eqref{e:4Admissible2},\eqref{e:4Admissible3},\eqref{e:4Admissible4}
from the proof of Example~\ref{x:Dyadic}ii).
Inductively construct a sequence of finite functions $F_n:\vec\Xi_{\mu'(n+1)}\to\THREE^{n}$ satisfying
\begin{gather}
\label{e:4Admissible5}
\xi'[\vec u\circ\Cantor] \;\subseteq\; \cball\big(\sbin\big(F_n(\vec u)\big),2^{-n-2}\big)  \\
\label{e:4Admissible6}
F_{n+1}(\vec v) \;\in\; F_{n}\big(\vec v|_{<\mu'(n+1)}\big) \circ\THREE, \quad \vec v\in\vec\Xi_{\mu'(n+2)} \enspace .
\end{gather}
meaning that $F_{n+1}$ `extends' $F_n$ by one signed digit $s_{n}$.
It follows that $F(\bar u)=\lim_n F_n\big(\bar u|_{<\mu'(n+1)}\big)$
satisfies $\xi'=\sbin\circ F$ and has modulus $\Lambda:m\mapsto\mu'(m+1)$.
Hence $G:=\tilde\sigma\circ F$ exhibits the claimed properties.
\\
Conveniently abbreviate with $F_0$ the constant empty string $\epsilon$.
For $n=1$ take $\vec v\in\vec\Xi_{\mu'(2)}$ and consider
$r(\vec v)\in\xi'[\vec v\circ\Cantor]\subseteq\cball\big(r(\vec v),2^{-3}\big)\subseteq[0;1]$
according to Equations~\eqref{e:4Admissible3} and \eqref{e:4Admissible4}.
Hence $|r(\vec v)-\sbin(F_0)|\leq \tfrac{3}{2}\cdot2^{-2}$; therefore 
\[ s_{0}\;:=\; \big\lfloor 2^{2}\cdot\big(r(\vec v)\:-\:\sbin(F_0)\big)\big\rceil\;\in\;\THREE \]
is well-defined; and $F_{1}(\vec v):=s_{0}$ has
$\sbin\big(F_{1}(\vec v)\big)=\tfrac{1}{2}+s_{0}\cdot2^{-2}$ with
$\big|r(\vec v)-\sbin\big(F_1(\vec v)\big)\big|\leq2^{-3}$, which asserts Equation~\eqref{e:4Admissible5}.
\\
Now turning to the induction step $F_{n+1}(\vec v)=F_n(\vec u)\circ s_{n}$ 
for $\vec u:=\vec v|_{<\mu'(n+1)}$ according to Equation~\eqref{e:4Admissible6}.
Note $r(\vec v)\in\xi'[\vec v\circ\Cantor]\subseteq\cball\big(r(\vec v),2^{-n\pmb{-3}}\big)$
and $\xi'[\vec v\circ\Cantor]\subseteq\xi'[\vec u\circ\Cantor]\subseteq\cball\big(r(\vec u),2^{-n\pmb{-2}}\big)$
according to Equations~\eqref{e:4Admissible2} and \eqref{e:4Admissible3} and \eqref{e:4Admissible4}.
Hence $|r(\vec v)-r(\vec u)|\leq2^{-n\pmb{-3}}$;
and $\big|r(\vec u)-\sbin(F_n)\big|\leq2^{-n-2}$ by induction hypothesis:
Yielding $\big|r(\vec u)\:-\:\sbin\big(F_n(\vec u)\big)\big|\leq\tfrac{3}{2}\cdot2^{-n\pmb{-2}}$.
Therefore
\[ s_{n}\;:=\;
\big\lfloor 2^{n+2}\cdot\big(r(\vec v)-\sbin\big(F_n(\vec u)\big)\big)\big\rceil\;\in\;\THREE \]
is well-defined; and $F_{n+1}(\vec v):=F_n(\vec u)\circ s_{n}$ has
$\sbin\big(F_{n+1}(\vec v)\big)=\sbin\big(F_{n}(\vec u)\big)+s_{n}\cdot2^{-n-2}$ with
$\big|r(\vec v)-\sbin(F_{n+1}(\vec v))\big|\leq2^{-n-3}$, 
which asserts Equation~\eqref{e:4Admissible5}
and completes the induction step.
\item[iii)]
To $r\in[0;1]$ consider any $\vec s=(s_0,\ldots,s_{n-1})\in\THREE^{n}$ with
$|r-\sbin(\vec s)|\leq2^{-n-1}$ according to (i').
This $\vec s$ extends to a signed binary expansion $\bar s$ of $r=\sbin(\bar s)$.
Similarly, the signed binary expansion of $(r'-r)\cdot2^n\in[-1/2;1/2]$ according to (i') 
yields $(s'_n,s'_{n+1},\ldots)\in\THREE^\IN$
with $r'-r=\sum\nolimits_{j\geq n} s'_j2^{-j-2}$.
Letting $s'_j:=s_j$ for $j<n$ results in $\bar s'$ with $\sbin(\bar s')=r'$.
Hence $D(\bar s,bar s')\leq2^{-n}$, and $\bar u:=\tilde\sigma(\bar s)$ and 
$\bar u':=\tilde\sigma(\bar s')$ have $D(\bar u,\bar u')\leq2^{-2n}$.
\item[iv)]
See Figure~\ref{f:Transducer}.
\qed\end{enumerate}\end{proof}
To summarize, designing a suitable encoding is challenging
already for the case $X=[0;1]$ of real numbers:
The immediate choice, binary expansion, turns out as computability-theoretical unsuitable, namely not admissible, recall Example~\ref{x:Binary}.
And even among the admissible ones, Example~\ref{x:Rational} is complexity-theoretically unsuitable.
Moreover, among those suitable for complexity considerations, 
some are only `polynomially admissible' (Example~\ref{x:Dyadic})
and some even `linearly admissible' (Example~\ref{x:SignDigit}):
as to be formalized in Definition~\ref{d:Admissible2} below.

\section{Quantitatively Admissible Representations on Cantor Space}
\label{s:Admissible2}

Guided by Examples~\ref{x:Dyadic} and \ref{x:SignDigit}, 
we now quantitatively refine qualitative Definition~\ref{d:Admissible}
polynomially and linearly, for the case of representations over Cantor space.
An initial idea might require such a representation to have polynomial/linear
modulus of continuity, instead of being quantitatively continuous;
however Lemma~\ref{l:Entropy}g) reveals that this condition is feasible
only for spaces of polynomial/linear entropy. Instead, Definition~\ref{d:Admissible2}
requires the representation to have modulus of continuity polynomial/linear `relative' 
to the entropy of the space under consideration:

\begin{definition}
\label{d:Admissible2}
For Landau's asymptotic $\calO(n)$,
abbreviate $\poly(n):=\calO(n)^{\calO(1)}$.
\begin{enumerate}
\item[a)]
Fix metric space $(X,d)$ with entropy $\Entropy$
as well as partial mappings $\xi,\xi':\subseteq U\to X$
with minimal moduli of continuity $\mu,\mu'$.
\item[a\,i)]
A \emph{polynomial reduction from $\xi'$ to $\xi$}
is a mapping $G:\dom(\xi')\to\dom(\xi)$
having modulus of continuity $\nu$
such that $\xi'\sqsubseteq\xi\circ G$
and $\nu\circ\mu\leq\mu'\circ\poly$.
\\
In this case we write ``$\xi'\reduceqP\xi$''.
\item[a\,ii)]
A \emph{linear reduction from $\xi'$ to $\xi$}
is a mapping $G:\dom(\xi')\to\dom(\xi)$
having modulus of continuity $\nu$
such that $\xi'\sqsubseteq\xi\circ G$
and $\nu\circ\mu\leq\mu'\circ\calO$.
\\
In this case we write ``$\xi'\reduceqL\xi$''.
\item[b)]
Call representation $\xi:\subseteq \Cantor\twoheadrightarrow X$ \emph{polynomially admissible} ~iff~ 
it satisfies (b\,i) and (b\,ii), where:
\item[b\,i)]
$\xi$ has modulus of continuity $\mu(n) \leq \Entropy\big(\poly(n)\big)$.
\item[b\,ii)] 
Every continuous surjective $\xi':\subseteq\Cantor\twoheadrightarrow X$ has $\xi'\reduceqP\xi$ from (a\,i).
\item[c)]
Call representation $\xi:\subseteq \Cantor\twoheadrightarrow X$ \emph{polynomially standard} ~iff~ 
it satisfies (b\,i) $\mu \leq \Entropy\circ\poly$ and 
(b\,ii') $\xi'\reduceqP\xi$ for every (not necessarily surjective) $\xi':\subseteq\Cantor\to X$
and the following condition (c\,i):
\item[c\,i)]
There exists a polynomial $p\in\IN[N]$ such that,
for every $n\in\IN$ and every $x,x'\in X$,
$d(x,x')\leq2^{-p(n)}$ implies $d_\xi(x,x')\leq2^{-\Entropy(n)}$
with the induced semi-metric $d_\xi$ according to Equation~\eqref{e:Modulus2}; \\
short: $d(x,x')\leq2^{-\poly(n)}\Rightarrow d_\xi(x,x')\leq2^{-\Entropy(n)}$.
\item[d)]
Call representation $\xi:\subseteq \Cantor\twoheadrightarrow X$ \emph{linearly admissible} ~iff~ 
it satisfies (d\,i) and (d\,ii), where:
\item[d\,i)]
$\xi$ has modulus of continuity $\mu(n) \leq \Entropy\big(\calO(n)\big)$.
\item[d\,ii)] 
Every continuous surjective $\xi':\subseteq\Cantor\twoheadrightarrow X$ has $\xi'\reduceqL\xi$ from (a\,ii).
\item[e)]
Call representation $\xi:\subseteq \Cantor\twoheadrightarrow X$ \emph{linearly standard} ~iff~ 
it satisfies (d\,i) $\mu \leq \Entropy\circ\calO$ and 
(d\,ii') $\xi'\reduceqL\xi$ for every (not necessarily surjective) $\xi':\subseteq\Cantor\to X$
and the following condition (e\,i):
\item[e\,i)]
$d(x,x')\leq2^{-\calO(n)} \quad\Rightarrow\quad d_\xi(x,x')\leq2^{-\Entropy(n)}$.
\end{enumerate}
\end{definition}
Definition~\ref{d:Admissible2}a) quantitatively refines qualitative continuous reduction from Definition~\ref{d:Admissible}a).
Again, any polynomially/linearly standard representation is \emph{a fortiori} polynomially/linearly admissible.
Conditions~(c) and (e) formalize that the semi-metric induced 
by the representation be `close' to the original metric of the space
in the sense of Example~\ref{x:PolyMetric}: as quantitative counterpart
to the qualitative \emph{final} topology in Definition~\ref{d:Admissible}c).

Observe that $\xi'\sqsubseteq\xi\circ G$ has modulus of continuity $\mu'\leq\nu\circ\mu$ by Remark~\ref{r:Modulus}c);
and Conditions~(a\,i) and (a\,ii) require this bound to be `almost' tight.
Similarly, Items~(b\,i) and (d\,i) quantitatively refine qualitative continuity from Definition~\ref{d:Admissible}b):
$X$ has entropy $\Entropy\leq\mu$ by Lemma~\ref{l:Entropy}a+g); and Items~(b\,i) and (d\,i) require this bound to be `almost' tight.
On both cases, `almost' tight is formalized such as to allow some (namely respectively linear and polynomial) slack in the argument, not in the value:
in order to make ``$\reduceqP$'' and ``$\reduceqL$'' partial orders and in particular transitive.

\begin{remark}
\label{r:Slack} 
Consider some inequality
\begin{equation}
\label{e:Inequality}
\forall n: \quad \mu(m)\;\leq\;\nu(n) \enspace . 
\end{equation}
between two moduli $\mu,\nu\in\IN^\IN$ according to Remark~\ref{r:Skolem}.
For example Lemma~\ref{l:Entropy}a+g) imposes such an inequality between the entropy of a space and the modulus of continuity of a representation over Cantor space;
and Remark~\ref{r:Modulus}c) imposes such an inequality on the modulus of continuity of a reduction from one representation to another.
Our general approach towards quantitatively refining qualitative conditions
takes such a necessary inequality, and then relaxes it:
from $\mu$ being perfectly bounded by $\nu$ as in the original Inequality~\eqref{e:Inequality}
to being \emph{linearly} or \emph{polynomially} bounded by $\nu$.
Example~\ref{x:PolyMetric} and Definitions~\ref{d:Entropy2} and \ref{d:Admissible3} below also follow these lines.

Inequality~\eqref{e:Inequality} can be weakened in various ways:
bounding $\mu$ by $\nu$ up to a constant offset only, additionally allow a constant factor (=linearly), or a constant power (polynomially).
These types of slack can be granted either on the values only, or on the arguments only, or both:
\begin{multline*}
\mu(m)\leq\nu(n)+\calO(1), \quad \mu(m)\leq\calO\big(\nu(n)\big), \quad \mu(m)\leq\poly\big(\nu(n)\big)  \\
\mu(m)\leq\nu\big(n+\calO(1)\big), \quad \mu(m)\leq\nu\big(\calO(n)\big), \quad \mu(m)\leq\nu\big(\poly(n)\big)  \\
\mu(m)\leq\nu\big(n+\calO(1)\big)+\calO(1), \quad \mu(m)\leq\calO\Big(\nu\big(\calO(n)\big)\Big), \quad \mu(m)\leq\poly\Big(\nu\big(\poly(n)\big)\Big)
\end{multline*}
Note that all these relations are transitive.
Restricted to strictly increasing mappings in $\IN^\IN$, 
the first line/type implies the second and coincides with the third.
Nevertheless, we make an effort to cover all non-decreasing (instead of only strictly increasing) moduli of 
continuity, such as in Equation~\eqref{e:LoInv}---until but excluding Section~\ref{s:Multifunc}.

Theorem~\ref{t:Standard} requires offsetting the argument.
\Martin{%
Moreover, quantitative reduction among representations (Definition~\ref{d:Admissible2}a) involves
composition with a third function, and that relation is transitive
only with the second type of slack: in the arguments.}
\end{remark}
Note the categorical similarity of qualitative/linear/polynomial admissibility 
to computable/polynomial-time \emph{completeness} in the discrete Theory of Computing:
all amount to being maximal with respect to a certain partial order.
Items~(c\,i) and (e\,i) refine Definition~\ref{d:Admissible}c\,i') by requiring
the semi-metric $d_\xi$ (rather than the topology) induced by $\xi$ to be polynomially/linearly related to $d$:

\begin{example}
\label{x:PolyMetric} 
Fix two metric spaces $(X,d)$ and $(X,d')$ with same domain $X$.
Recall that $d$ and $d'$ are called \emph{equivalent} if they induce the same topology;
and \emph{strongly} equivalent if they satisfy both $d\leq \calO(d')$ and $d'\leq\calO(d)$.
\\
For bounded and for \Martin{complete} convex spaces $X$, 
strong equivalence can be rephrased as follows:
\begin{enumerate}
\item[i)] $\displaystyle 
\forall n\in\IN \; \forall x,x'\in X: \quad
d(x,x')\leq 2^{-(n+\calO(1))} \;\Rightarrow\; d'(x,x')\leq2^{-n} 
\quad\wedge\quad
d'(x,x')\leq2^{-(n+\calO(1))} \;\Rightarrow\; d(x,x')\leq2^{-n}$.
\end{enumerate}
This suggests the following notions `interpolating' between equivalence and strong equivalence:
\begin{enumerate}
\item[ii)] Call $d$ and $d$' \emph{linearly equivalent} if
$d(x,x')\leq 2^{-\calO(n)}$ implies $d'(x,x')\leq2^{-n}$
and $d'(x,x')\leq2^{-\calO(n)}$ implies $d(x,x')\leq2^{-n}$.
\item[iii)] Call $d$ and $d$' \emph{polynomially equivalent} if 
$d(x,x')\leq 2^{-\poly(n)}$ implies $d'(x,x')\leq2^{-n}$
and $d'(x,x')\leq2^{-\poly(n)}$ implies $d(x,x')\leq2^{-n}$.
\end{enumerate}
More generally, consider metric $d'$ `rescaled' by some nondecreasing unbounded $\nu:\IN\to\IN$ similarly to Lemma~\ref{l:Entropy}j),
and $d$ rescaled by some $\mu$, and compare these in the sense of (i),(ii),(iii) as follows:
\begin{enumerate}
\item[i')]
Call $d$ \emph{strongly $(\mu,\nu)$-equivalent} to $d'$ if 
$d(x,x')\leq 2^{-\mu(n+\calO(1))}$ implies $d'(x,x')\leq2^{-\nu(n)}$ and if also
$d'(x,x')\leq2^{-\nu(n+\calO(1))}$ implies $d(x,x')\leq2^{-\mu(n)}$.
\item[ii')]
Call $d$ \emph{linearly $(\mu,\nu)$-equivalent} to $d'$ if 
$d(x,x')\leq 2^{-\mu(\calO(n))}$ implies $d'(x,x')\leq2^{-\nu(n)}$ and if also
$d'(x,x')\leq2^{-\nu(\calO(n))}$ implies $d(x,x')\leq2^{-\mu(n)}$.
\item[iii')]
Call $d$ \emph{polynomially $(\mu,\nu)$-equivalent} to $d'$ if 
$d(x,x')\leq 2^{-\mu(\poly(n))}$ implies $d'(x,x')\leq2^{-\nu(n)}$ and if also
$d'(x,x')\leq2^{-\nu(\poly(n))}$ implies $d(x,x')\leq2^{-\mu(n)}$.
\end{enumerate}
In case $\mu=\id=\nu$, (i')--(iii') boil down to (i)--(iii).
These relations\footnote{%
We consider polynomial `slack' in the argument $\mu(\poly(n))$, 
not of the value $\poly(\mu(n))$: Making the latter transitive involves
more subtle quantification over polynomials in the antecedent and consequent,
and introduces third $\poly$'s in Theorem~\ref{t:Main1}(d\,i) and (e\,i);
similarly for the linear case.} 
are `generalized' symmetric and transitive:
If $d$ is strongly/linearly/polynomially $(\mu,\nu)$-equivalent to $d'$
then $d'$ is strongly/linearly/polynomially $(\nu,\mu)$-equivalent to $d$;
and if additionally $d''$ is strongly/linearly/polynomially $(\nu,\lambda)$-equivalent to $d'$,
then it is strongly/linearly/polynomially $(\mu,\lambda)$-equivalent to $d$. 
\\
In view of Equation~\ref{e:Modulus2},
one can thus concisely reformulate Definition~\ref{d:Admissible2}(b\,i) and (c\,i)  
as polynomial $(\id,\Entropy)$-equivalence of $d$ and $d_\xi$; 
and (d\,i)+(e\,i) mean linear $(\id,\Entropy)$-equivalence of $d$ and $d_\xi$.
\end{example}
Definition~\ref{d:Admissible2}c\,i) can alternatively be rephrased as:
Any $x,x'$ with $d(x,x')\leq2^{-\poly(n)}$ have respective $\xi$-names $\bar u,\bar u'$
with $D(\bar u,\bar u')\leq2^{-\Entropy(n)}$. 
According to Example~\ref{x:Dyadic}, 
the dyadic representation $\delta$ of $[0;1]$
is therefore polynomially standard;
and according to Example~\ref{x:SignDigit},
the signed-digit representation $\sigma$ of $[0;1]$
is linearly standard.

\subsection[Quantitative \emph{Main Theorem}]{Quantitative \emph{Main Theorem} for Cantor Space Representations}
\label{ss:Main1}

We can now state a first quantitative variant of the qualitative Kreitz-Weihrauch \emph{Main Theorem} 
(Fact~\ref{f:KreitzWeihrauch}), here for Cantor space representations: 

\begin{theorem}
\label{t:Main1}
\begin{enumerate}
\item[a\,i)]
Suppose $(X,d)$ admits a polynomially standard representation $\xi:\subseteq\Cantor\twoheadrightarrow X$
and $\xi':\subseteq\Cantor\twoheadrightarrow X$ is polynomially admissible.
Then $\xi'$ is polynomially standard, too.
\item[a\,ii)]
Suppose $(X,d)$ admits a linearly standard representation $\xi:\subseteq\Cantor\twoheadrightarrow X$
and $\xi':\subseteq\Cantor\twoheadrightarrow X$ is linearly admissible.
Then $\xi'$ is linearly standard, too.
\item[b\,i)]
Fix a polynomially standard representation $\xi:\subseteq\Cantor\twoheadrightarrow X$ 
of metric space $(X,d)$ with entropy $\Entropy$. Consider $f:X\to Y$ for metric space $(Y,e)$.
If $f\circ\xi$ has modulus of continuity $\Entropy\circ\lambda$, then $f$ has modulus of continuity $\poly\circ\lambda$.
\item[b\,ii)]
If $\xi:\subseteq\Cantor\twoheadrightarrow X$ is linearly standard,
then $f$ has modulus of continuity $\calO\circ\lambda$.
\item[c)]
Let $\xi:\subseteq\Cantor\twoheadrightarrow X$ and $\upsilon:\subseteq\Cantor\twoheadrightarrow Y$
denote representations of metric spaces $(X,d)$ and $(Y,e)$ 
with respective entropies $\Entropy$ and $\theta$.
\item[c\,i)]
Suppose $\xi$ and $\upsilon$ are polynomially standard.
Consider $f:X\to Y$ with modulus of continuity $\lambda$.
Then $f$ admits a $(\xi,\upsilon)$-realizer $F:\dom(\xi)\to\dom(\upsilon)$
with modulus of continuity $\Lambda$ such that 
$\Lambda\circ\theta\leq \Entropy\circ\poly\circ\lambda\circ\poly$.
\item[c\,ii)]
Suppose $\xi$ and $\upsilon$ are linearly standard.
Consider $f:X\to Y$ with modulus of continuity $\lambda$.
Then $f$ admits a $(\xi,\upsilon)$-realizer $F$
with modulus of continuity $\Lambda$ such that 
$\Lambda\circ\theta\leq \Entropy\circ\calO\circ\lambda\circ\calO$.
\item[d)]
Let $\xi:\subseteq\Cantor\twoheadrightarrow X$ and $\upsilon:\subseteq\Cantor\twoheadrightarrow Y$
denote representations of metric spaces $(X,d)$ and $(Y,e)$ 
with respective entropies $\Entropy$ and $\theta$.
\item[d\,i)]
Suppose $\xi$ and $\upsilon$ are polynomially standard.
Let $f:X\to Y$ admit a $(\xi,\upsilon)$-realizer $F$ with modulus of continuity $\Entropy\circ\lambda$.
Then $f$ has modulus of continuity $\poly\circ\lambda\circ\theta\circ\poly$.
\item[d\,ii)]
Suppose $\xi$ and $\upsilon$ are linearly standard.
Let $f:X\to Y$ admit a $(\xi,\upsilon)$-realizer $F$ with modulus of continuity $\Entropy\circ\lambda$.
Then $f$ has modulus of continuity $\calO\circ\lambda\circ\theta\circ\calO$.
\end{enumerate}
\end{theorem}
The proof is deferred to Subsection~\ref{ss:Proofs1}.
Fact~\ref{f:KreitzWeihrauch} relates qualitative continuity of functions to qualitative continuity of their realizers;
and Theorem~\ref{t:Main1}c+d) does so quantitatively---relative to the entropies of their co/domains:
and up to polynomial/linear `slack' in both the argument \emph{and value}.

\begin{corollary}
\label{c:Main1}
Fix some non-decreasing unbounded $\lambda:\IN\to\IN$.
Since $[0;1]$ has linear entropy, 
Theorem~\ref{t:Main1}c+d) together with Examples~\ref{x:Dyadic} and \ref{x:SignDigit} implies:
\begin{enumerate}
\item[i)]
A function $f:[0;1]\to[0;1]$ has modulus of continuity $\poly\circ\lambda\circ\poly$
~iff~ it has a $(\delta,\delta)$-realizer with modulus of continuity $\poly\circ\lambda\circ\poly$.
\item[ii)]
A function $f:[0;1]\to[0;1]$ has modulus of continuity $\calO\circ\lambda\circ\calO$
~iff~ it has a $(\sigma,\sigma)$-realizer with modulus of continuity $\calO\circ\lambda\circ\calO$.
\end{enumerate}
\end{corollary}
This improves \cite[Theorem~2.19]{Ko91}; see also \cite[Theorem~14]{DBLP:conf/lics/KawamuraS016}.

\subsection{Relative Entropy, Capacity, and (Almost) Homogeneous Spaces}
\label{ss:Entropy2}

Corollary~\ref{c:Main1}ii) is tailored to the reals with the linearly admissible representation $\sigma$ over Cantor as an ultrametric space.
Note that for both the reals and Cantor space, the entropy of a ball/interval depends only on its radius but not on its center.
Definition~\ref{d:Entropy2}d) calls this property \emph{perfect homogenity}, and then proceeds to weakens it:
towards generalizing Corollary~\ref{c:Main1}ii) to larger classes of both represented and ground spaces. 

\begin{definition}
\label{d:Entropy2}
Fix a compact metric space $(X,d)$ with unary and binary entropy 
$\unary{\Entropy}_X$ and $\Entropy_X=\lceil\log_2\unary{\Entropy}_X\rceil$, respectively.
\begin{enumerate}
\item[a)]
For $m\in\IN$ abbreviate with $\unary{\Entropy}_{m}=\sup_{x\in X} \unary{\Entropy}_{\cball_m(x)}$
the (maximum) \emph{relative} unary entropy of all closed balls of radius $2^{-m}$;
write $\Entropy_m=\lceil\log_2\unary{\Entropy}_{m}\rceil$ for their (maximum) relative binary entropy.
\item[b)]
The unary \emph{capacity} of $X$ is
$\unary{\Capacity}=\unary{\Capacity}_X:\IN\to\IN$ such that $\unary{\Capacity}(n)$, but no more,
points of pairwise distance $>2^{-n}$ fit into $X$. \\
Write $\Capacity_X=\lfloor\log_2\unary{\Capacity}_X\rfloor$ for the (binary) capacity.
\item[c)]
For $m\in\IN$ abbreviate with $\unary{\Capacity}_{m}=\inf_{x\in X} \unary{\Capacity}_{\cball_m(x)}$
the (minimum) \emph{relative} unary entropy of all closed balls of radius $2^{-m}$, considered as subspaces of $X$; 
write $\Capacity_{m}=\lceil\log_2\unary{\Capacity}_{m}\rceil$ for the (minimum) relative binary capacity.
\item[d)]
Call $(X,d)$ \emph{(perfectly) homogeneous} if it holds both $\unary{\Entropy}_{\cball_m(x)}=\unary{\Entropy}_m$ 
and $\unary{\Capacity}_{\cball_m(x)}=\unary{\Capacity}_m$ for all $x\in X$ and all $m\in\IN$.
\item[d')]
$X$ is \emph{strongly homogeneous} if it holds $\forall x\in X: \unary{\Entropy}_{m}(n) \leq\unary{\Entropy}_{\cball_m(x)}\big(n+\calO(1)\big)$
with the constant in $\calO$ independent of $x\in X$ and of $n,m\in\IN$.
\\
$X$ is \emph{linearly homogeneous} if it holds $\forall x\in X: \unary{\Entropy}_{m}(n)\leq\unary{\Entropy}_{\cball_m(x)}\big(\calO(n)\big)$.
\\
$X$ is \emph{polynomially homogeneous} if it holds $\forall x\in X: \unary{\Entropy}_{m}(n)\leq\unary{\Entropy}_{\cball_m(x)}\big(\poly(n)\big)$.
\end{enumerate}
\end{definition}
Entropy refers to minimal coverings, capacity to maximal packings; cmp. \cite[Definition~6.2]{Wei03}.
Accordingly we round the binary entropy up and round the binary capacity down to the next integer power of two.
All compact metric groups are perfectly homogeneous \cite{DBLP:conf/csl/PaulyS020}.
Item~(d') about polynomial/linear/strong (but not perfect) 
homogeneity can equivalently be formulated in terms of the unary (but not the binary) capacity instead of the unary (but not the binary) entropy,
according to Item~(a) in the following Lemma: 

\begin{lemma}
\label{l:Entropy2}
\begin{enumerate}
\item[a)]
It holds $\unary{\Entropy}(n)\leq\unary{\Capacity}(n)\leq\unary{\Entropy}(n+1)$; \\ 
hence $\Entropy(n)\leq\Capacity(n)+1\leq\Entropy(n+1)+1$.
\item[b)] 
If $(X,d)$ and $(Y,e)$ have unary entropies $\unary{\Entropy}_X$ and $\unary{\Entropy}_Y$,
then $(X\times Y,\max\{d,e\})$ has unary entropy $\unary{\Entropy}_{X\times Y}\leq \unary{\Entropy}_X\cdot \unary{\Entropy}_Y$,
possibly $<$.
\\
If $(X,d)$ and $(Y,e)$ have unary capacities $\unary{\Capacity}_X$ and $\unary{\Capacity}_Y$,
then $(X\times Y,\max\{d,e\})$ has unary capacity $\unary{\Capacity}_{X\times Y}\geq \unary{\Capacity}_X\cdot \unary{\Entropy}_Y$,
possibly $>$.
\\
If both $(X,d)$ and $(Y,e)$ are polynomially/linearly/strongly homogeneous, then so is $(X\times Y,\max\{d,e\})$. 
\item[c)] 
Let compact $(X_j,d_j)$ have diameter $\leq1$ and unary entropy $\unary{\Entropy}_j$ 
and unary capacity $\unary{\Capacity}_j$ for each $j\in\IN$.
Then $\big(\prod_jX_j,\sup_j d_j/2^j\big)$ is compact and has unary entropy 
$\unary{\Entropy}(n)\leq\prod_{j\leq n} \unary{\Entropy}_j(n-j)$
and unary capacity $\unary{\Capacity}(n)\geq\prod_{j\leq n}\unary{\Capacity}_j(n-j)$. \\
Hence the binary entropy satisfies ${\Entropy}(n)\leq\sum_{j\leq n} {\Entropy}_j(n-j)$, 
and the binary capacity satisfies ${\Capacity}(n)\geq\sum_{j\leq n}{\Capacity}_j(n-j)$. 
\item[d)]
Suppose all $(X_j,d_j)$ are \emph{uniformly} polynomially/linearly/strongly homogeneous
in the sense that there exists some $c\in\IN$ satisfying, for all $j,n,m\in\IN$ and $x\in X_j$, 
\[ \unary{\Entropy}_{j,m}(n) \quad\leq\quad 
\unary{\Entropy}_{j,\cball_m(x)}(n+c),  \qquad
\leq \unary{\Entropy}_{j,\cball_m(x)}(c+c\cdot n),  \qquad
\leq \unary{\Entropy}_{j,\cball_m(x)}(c+c\cdot n^c) \enspace , \]
respectively.
Then $\big(\prod_jX_j,\sup_j d_j/2^j\big)$ is again polynomially/linearly/strongly homogeneous.
\item[e)]
Every compact metric space $(X,d)$ satisfies $\unary{\Entropy}_X(n)\leq\unary{\Entropy}_X(m)\cdot\unary{\Entropy}_{m}(n)$
for all (!) $n,m\in\IN$. Hence $\Entropy_X(n)\leq{\Entropy}_X(m)+{\Entropy}_{m}(n)$.  
\item[f)]
If $(X,d)$ is perfectly homogeneous ultrametric, then the inequality from (e) 
becomes equality: $\unary{\Entropy}_X(n)=\unary{\Entropy}_X(m)\cdot\unary{\Entropy}_{m}(n)$.
And rounding up to integer powers of two yields $\Entropy_X(n)+1\geq{\Entropy}_X(m)+{\Entropy}_{m}(n)$.
\item[f')]
If $(X,d)$ is strongly homogeneous ultrametric, then it holds
$\unary{\Entropy}_X\big(n+\calO(1)\big)\geq\unary{\Entropy}_X(m)\cdot\unary{\Entropy}_{m}(n)$.
And ${\Entropy}_X\big(n+\calO(1)\big)+1\geq{\Entropy}_X(m)+{\Entropy}_{m}(n)$.
\\
If $(X,d)$ is linearly homogeneous ultrametric, then 
$\unary{\Entropy}_X\big(\calO(n)\big)\geq\unary{\Entropy}_X(m)\cdot\unary{\Entropy}_{m}(n)$.
And ${\Entropy}_X\big(\calO(n)\big)+1\geq{\Entropy}_X(m)+{\Entropy}_{m}(n)$.
\\
If $(X,d)$ is polynomially homogeneous ultrametric, then 
$\unary{\Entropy}_X\big(\poly(n)\big)\geq\unary{\Entropy}_X(m)\cdot\unary{\Entropy}_{m}(n)$.
And ${\Entropy}_X\big(\poly(n)\big)+1\geq{\Entropy}_X(m)+{\Entropy}_{m}(n)$.
\item[g)]
Every compact ultrametric space $(X,d)$ satisfies $\unary{\Capacity}_X(n)\geq\unary{\Capacity}_X(m)\cdot\unary{\Capacity}_{m}(n)$; 
and hence $\Capacity_X(n)\geq{\Capacity}_X(m)+{\Capacity}_{m}(n)$ for all $n>m$.
\item[h)]
If $(X,d)$ is perfectly homogeneous ultrametric, then the inequality from (g) becomes 
an equality: $\unary{\Capacity}_X(n)=\unary{\Capacity}_X(m)\cdot\unary{\Capacity}_{m}(n)$ for all $n>m$.
And rounding down to integer powers of two yields $\Capacity_X(n)\leq{\Capacity}_X(m)+{\Capacity}_{m}(n)+1$. 
\item[h')] 
If $(X,d)$ is strongly homogeneous ultrametric, then 
$\unary{\Capacity}_X(n)\leq\unary{\Capacity}_X(m)\cdot\unary{\Capacity}_{m}\big(n+\calO(1)\big)$; 
and hence $\Capacity_X(n)\leq{\Capacity}_X(m)+{\Capacity}_{m}\big(n+\calO(1)\big)+1$ for all $n>m$.
\\
If $(X,d)$ is linearly homogeneous ultrametric, then 
$\unary{\Capacity}_X(n)\leq\unary{\Capacity}_X(m)\cdot\unary{\Capacity}_{m}\big(\calO(n)\big)$;
and hence $\Capacity_X(n)\leq{\Capacity}_X(m)+{\Capacity}_{m}\big(\calO(n)\big)+1$ for all $n>m$.
\\
If $(X,d)$ is polynomially homogeneous ultrametric, then 
$\unary{\Capacity}_X(n)\leq\unary{\Capacity}_X(m)\cdot\unary{\Capacity}_{m}\big(\poly(n)\big)$;
and hence $\Capacity_X(n)\leq{\Capacity}_X(m)+{\Capacity}_{m}\big(\poly(n)\big)+1$ for all $n>m$.
\end{enumerate}
\end{lemma}
Lemma~\ref{l:Entropy}d) now follows from Items~(a)+(d).
Note that Claim~(a) only refers to the (absolute) entropy/capacity.
To see the first part of (a), take a maximal set $X_n\subseteq X$ of points of pairwise distance $>2^{-n}$.
Then the closed balls $\cball_n$ with centers $x\in X_n$ cover $X$: 
any point $x'\in X\setminus\bigcup_{x\in X_n} \cball(x,2^{-n})$ would contradict the maximality of $X_n$.
For the second part of (a) observe that any $\cball_{n+1}$ can cover at most one of the points in $X_n$.
Finally recall that the binary capacity rounds down while the binary entropy rounds up.
\\
To see (e), abbreviate $J=\unary{\Entropy}(m)$ and let $x_1,\ldots,x_J$ denote centers of an optimal covering of $X$ by $\cball_m$s.
Then, $\cball_n$s covering the $\cball_m(x_j)$, together cover $\bigcup_j\cball_m(x_j)\supseteq X$; hence
\begin{equation}
\label{e:Entropy2}
\unary{\Entropy}_X(n) \;\leq\; \sum\nolimits_{j<J} \unary{\Entropy}_{\cball_m(x_j)}(n) 
\;\leq\; J\cdot\unary{\Entropy}_{m}(n) \enspace . 
\end{equation}
For ultrametric spaces, the first inequality in Equation~\eqref{e:Entropy2} becomes an equality:
because then minimal coverings are unique, see Observation~\ref{o:Ultra}v) below.
And the second inequality in Equation~\eqref{e:Entropy2} becomes an equality for perfectly homogeneous spaces;
hence Item~(f) follows, and (f') similarly.
\\
Observation~\ref{o:Ultra}iv) shows that every compact ultrametric space satisfies
\begin{equation}
\label{e:Capacity2}
\unary{\Capacity}_X(n) \;=\; \sum\nolimits_{j<J} \unary{\Capacity}_{\cball(x_j,2^{-m})}(n) 
\;\geq\; J\cdot \unary{\Capacity}_m(n)  \enspace ,
\end{equation}
where $J=\unary{\Capacity}_X(m)$ and $x_1,\ldots,x_J$ denote a maximal packing in $X$ with pairwise distances $>2^{-m}$.
This implies (g); and the second inequality in Equation~\eqref{e:Capacity2} becomes an equality 
for perfectly homogeneous spaces, asserting (h), and (h').
The above arguments have employed some of the following common properties of (compact) ultrametric spaces \cite[\S2+\S5]{Ultrametric}:

\begin{observation}
\label{o:Ultra}
Let $(\calU,D)$ denote a compact ultrametric space.
Then every non-empty closed subset of $\calU$ is again a compact ultrametric space.
The closed balls $\cball(x,r)$ of radius $r>0$, $x\in\calU$, 
\begin{enumerate}
\item[i)] are topologically open, and finitely many of them cover $\calU$.
\item[ii)] If $D(x,x')\leq r$, then the two balls $\cball(x,r),\cball(x',r)$ are equal;
\item[iii)] If $D(x,x')>r$, then the two balls $\cball(x,r),\cball(x',r)$ are disjoint:
\item[iv)] In the latter case, all $z\in\cball(x,r)$ and $z'\in\cball(x',r)$ satisfy $D(z,z')=D(x,x')$.
\item[v)] Let $\calB_r$ and $\calB'_r$ denote two (not necessarily finite) sets of closed balls 
  $\cball(x,r)$ of radius $r>0$, $x\in\calU$. Suppose their unions agree: $\bigcup\calB'_r=\bigcup\calB_r$.
  Then it holds $\bigcup\calB_r=\bigcup(\calB_r\cap\calB'_r)$. 
\end{enumerate}
In particular \emph{minimal} covers by balls of same radius are unique
(but their centers are of course not).
\end{observation}
Recall that, for a set $\calY$ of subsets $Y\subseteq X$, $\bigcup\calY$ means $\bigcup_{Y\in\calY} Y$.

\begin{example}
\label{x:Entropy2}
\begin{enumerate}
\item[a)]
In Cantor space, any two balls of same radius $r>0$ are isometric: perfect homogeneity.
The relative unary entropy/capacity here is $\unary{\Entropy}_m(n)=\max\{2^{n-m},1\}=\unary{\Capacity}_m(n+1)$,
relative binary entropy/capacity ${\Entropy}_m(n)=\max\{n-m,0\}={\Capacity}_m(n+1)$.
Hence most inequalities in Lemma~\ref{l:Entropy2} become equalities:
$\unary{\Entropy}(n)=\unary{\Entropy}(m)\cdot\unary{\Entropy}_{m}(n)$,
${\Entropy}(n)={\Entropy}(m)+{\Entropy}_{m}(n)$, and
$\unary{\Capacity}(n)=\unary{\Capacity}(m)\cdot\unary{\Capacity}_{m}(n)$,
${\Capacity}(n)={\Capacity}(m)+{\Capacity}_{m}(n)$ for $n>m$.
\item[b)] 
In the space $\dual{(\TWO^*)}$ of total string predicates (Example~\ref{x:Example}f), too,
all balls of same radius $r>0$ are isometric: perfect homogeneity.
The relative binary entropy here is $\Entropy_m(n)=\max\{2^n-2^m,0\}=\Capacity_m(n+1)$; cmp. Lemma~\ref{l:Entropy}b).
Again, most inequalities in Lemma~\ref{l:Entropy2} become equalities, such as
${\Entropy}(n)={\Entropy}(m)+{\Entropy}_{m}(n)$ and ${\Capacity}(n)={\Capacity}(m)+{\Capacity}_{m}(n)$ for $n>m$.
\item[c)]
The real unit interval $[0;1]$ with unary entropy $\unary{\Entropy}(n)=\max\{1,2^{n-1}\}$
and unary capacity $\unary{\Capacity}(n)=2^n$ is strongly homogeneous: 
Different balls $\cball(x,r)$ of same radius $r=1/2^m$ need not be isometric, since those close to the the boundary 
of $[0;1]$ may get `cut off'---but by at most half, 
whose effect on the entropy is bounded by a constant shift $n\mapsto n+1$:
$\unary{\Entropy}_{\cball(x,1/2^m)}(n)\leq \unary{\Entropy}_m(n+1)$
and similarly $\unary{\Capacity}_{\cball(x,1/2^m)}(n)\leq \unary{\Capacity}_{m}(n+1)$
for $n>m$.
\\
Moreover the relative binary entropy and capacity satisfy
$\Entropy_m(n)=n-m=\Capacity_m(n)$ for $n>m>0$; 
hence ${\Entropy}(n)={\Entropy}(m)+{\Entropy}_{m}(n)$
and $\Capacity(n)=\Capacity(m)+\Capacity_m(n)$.
\end{enumerate}
\end{example}

\begin{remark}
\label{r:OpenClosed}
\begin{enumerate}
\item[a)]
Definition~\ref{d:Entropy2}a) about the relative (say, unary) entropy $\unary{\Entropy}_{X'}$, $X':=\cball_m(x)$, is actually ambiguous:
It might refer to considering $X'$ as a metric subspace of $X$, and thus to be covered by balls with centers in $X'$;
while a more relaxed conception could allow $X'$ to be covered by (possibly fewer) balls with centers from entire $X$.
Note that the relative capacity does not suffer from such ambiguity.
In fact all considerations in the present work apply to both perspectives of relative entropy.
In particular Lemma~\ref{l:Entropy2}a) shows that the two variants cannot be too far from each other.
\item[b)]
Definition~\ref{d:Entropy2} considers entropy
in terms of \emph{closed} balls of radius $1/2^n$,
and capacity referring to distances \emph{strictly} greater than $1/2^n$.
Lemma~\ref{l:Entropy2}a) holds also for the different combination,
with entropy defined in terms of \emph{open} balls and capacity 
as distances greater that \emph{or equal} to $1/2^n$.
Again, the two notions differ by at most one in argument/value.
\item[c)]
Alternatively to Equation~\eqref{e:Modulus}, one could consider a variant modulus of continuity
such as to satisfy \emph{strict} inequalities in both hypothesis and conclusion.
Again, the two notions differ by at most one in argument/value.
Mixing non-/strict inequalities in hypothesis and conclusion
on the other hand would violate the composition rule in Remark~\ref{r:Modulus}c).
\end{enumerate}\end{remark}

\subsection{Constructing Quantitatively Standard Representations}
\label{ss:Standard}

Theorem~\ref{t:Main1} supposes that the space(s) under consideration
be equipped with linearly/polynomially standard representations to begin with.
We show that many spaces indeed do admit such representations.

\begin{theorem}
\label{t:Standard} 
Fix a compact metric space $(X,d)$ with binary entropy $\Entropy$
and a strictly increasing positive integer sequence $n_k$, $k\in\IN$.
Abbreviate $n_{-1}:=-\infty$ and $\mu(n):=\Entropy(1+n_0)+\sum_{k=0}^K \Entropy_{n_{k}}\big(1+n_{k+1}\big)$ for $n_{K-1}<n\leq n_K$.
\begin{enumerate}
\item[i)]
Suppose it holds $n_{k+1}\leq\calO(n_k)$ and $\mu(n_K)\leq \Entropy_X\big(\calO(n_K)\big)$.
Then $X$ admits a linearly standard representation $\xi:\subseteq\Cantor\twoheadrightarrow X$ with modulus of continuity $\mu$.
\item[ii)]
Suppose it holds $n_{k+1}\leq\poly(n_k)$ and
$\mu(n_K) \leq \Entropy_X\big(\poly(n_K)\big)$.
Then $X$ admits a polynomially standard representation $\xi:\subseteq\Cantor\twoheadrightarrow X$ with modulus of continuity $\mu$.
\end{enumerate}
\end{theorem}
Consider for instance the real unit interval $[0;1]$ with $n_k:\equiv k$,
hence $\mu(K)=\Entropy(1)+\sum_{k=0}^K \Entropy_{k}(k+2)=0+2\cdot(K+1)$ according to Example~\ref{x:Entropy2}: 
Item~(i) here recovers, and thus generalizes, the (properties of the) signed binary representation from Example~\ref{x:SignDigit}.
Item~(ii) recovers and generalizes the following statement implicit in \cite[\S3.1]{DBLP:conf/lics/KawamuraS016}.

\begin{corollary}
\label{c:Standard}
Let $(X,d)$ denote a compact metric space whose (non-relative binary) entropy grows at least and at most polynomially:
\[ \exists C>0: \quad n^{1/C}/C \;\leq\; \Entropy(n) \;\leq\; C\cdot n^C \enspace .\]
Then $X$ admits a polynomially standard representation.
\end{corollary}
Indeed, again take $n_k:\equiv k$ and observe 
\begin{multline*}
\sum\nolimits_{k=0}^K \Entropy_k(k+2)\;\leq\; \sum\nolimits_{k=0}^K \Entropy(k+2) \;\leq\; 
(K+1)\cdot\Entropy(K+2) \;\leq \\
\leq\; C\cdot(K+2)^{1+C}  \;=\; n^{1/C}/C \;\leq\; \Entropy(n) \;\leq\; \Entropy\big(\poly(K)\big) 
\end{multline*}
for $N:=C\cdot(K+2)^{1+C}$ and $n:=(C\cdot N)^C\leq\poly(K)$.

\medskip
By Lemma~\ref{l:Entropy}h),
every connected metric space has at least linear binary entropy;
however that lower bound is not sufficient for Corollary~\ref{c:Standard}:
Consider the generalized Hilbert cube
$\prod_n [0; 2^{-a_n}]$ equipped with the supremum norm,
where $(a_n) \subseteq \IN$ is non-decreasing and unbounded.
If $(a_n)$ grows slowly, $\prod_n [0; 2^{-a_n}]$ will have `large' entropy.
Taking $(a_n)$ in a careful way yields a connected space which
violates the prerequisites of Theorem~\ref{t:Standard}.

On the other hand the class of spaces admitting linearly/polynomially standard representations is closed under Cartesian products:

\begin{theorem}
\label{t:Cartesian}
\begin{enumerate}
\item[a)]
Fix compact metric spaces $(X,d)$ and $(Y,e)$
with representations 
$\xi:\subseteq\Cantor\twoheadrightarrow X$ and
$\upsilon:\subseteq\Cantor\twoheadrightarrow Y$.
These induce a representation $\xi\times\upsilon:\subseteq\Cantor\twoheadrightarrow X\times Y$
of $\big(X\times Y,\max\{d,e\}\big)$ with the following property:
\item[i)]
If both $\xi$ and $\upsilon$ are linearly admissible, then so is $\xi\times\upsilon$.
\item[ii)]
If both $\xi$ and $\upsilon$ are linearly standard, then so is $\xi\times\upsilon$.
\item[iii)]
If both $\xi$ and $\upsilon$ are polynomially admissible, then so is $\xi\times\upsilon$.
\item[iv)]
If both $\xi$ and $\upsilon$ are polynomially standard, then so is $\xi\times\upsilon$.
\item[b)]
Fix compact metric spaces $(X_j,d_j)$ of entropies $\Entropy_j$ \Martin{and diameters $1/2^j\leq\diam(X_j)\leq1$}, 
with representations $\xi_j:\subseteq\Cantor\twoheadrightarrow X_j$, $j\in\IN$. 
These induce a representation $\prod_j\xi_j:\subseteq\Cantor\twoheadrightarrow \prod_j X_j$
of $\big(\prod_j X_j,\max d_j/2^n\big)$ with the following property:
\item[b\,i)]
Suppose $\xi_j$ are \emph{uniformly} linearly admissible:
in that there exists some $c\in\IN$ independent of $j$ such that
$\xi_j$ has modulus of continuity $\mu_j(n)\leq \Entropy_j(c+c\cdot n)$;
and to every representation $\xi'_j:\subseteq\Cantor\twoheadrightarrow X_j$
with modulus of continuity $\mu'_j$
there exists a mapping $F_j:\dom(\xi'_j)\to\dom(\xi_j)$
with modulus of continuity $\lambda_j$ such that
$\xi'_j=\xi_j\circ F_j$ and 
$\lambda_j\circ\mu_j(n)\leq\mu_j'(c+c\cdot n)$.
Then $\prod_j\xi_j$ is again linearly admissible.
\item[b\,ii)]
Suppose $\xi_j$ are \emph{uniformly} linearly standard in the following sense:
\\
To every (not necessarily surjective) family of mappings $\xi'_j:\subseteq\Cantor\to X_j$
with modulus of continuity $\mu'_j$
there exists a mapping $F_j:\dom(\xi'_j)\to\dom(\xi_j)$
with modulus of continuity $\lambda_j$ such that
$\xi'_j=\xi_j\circ F_j$ and $\lambda_j\circ\mu_j\leq\mu_j'\circ\calO$.
\\
And, to any $(x_j),(x'_j)\in\prod_j X_j$ with $d_j(x_j,x'_j)\leq2^{-c-c(n-j)}$,
there exist $\bar u_j,\bar u'_j\in\Cantor$ with $\xi_j(\bar u_j)=x_j$
and $\xi_j(\bar u'_j)=x'_j$ and $D(\bar u,\bar u')\leq2^{-\Entropy_j(n)}$.
\\
Then $\prod_j\xi_j$ is again linearly standard.
\item[iii)]
Suppose $\xi_j$ are \emph{uniformly} polynomially admissible:
in that there exists some $c\in\IN$ independent of $j$ such that
$\xi_j$ has modulus of continuity $\mu_j(n)\leq \Entropy_j(c+c\cdot n^c)$;
and to every representation $\xi'_j:\subseteq\Cantor\twoheadrightarrow X_j$
with modulus of continuity $\mu'_j$
there exists a mapping $F_j:\dom(\xi'_j)\to\dom(\xi_j)$
with modulus of continuity $\lambda_j$ such that
$\xi'_j=\xi_j\circ F_j$ and 
$\lambda_j\circ\mu_j(n)\leq\mu_j'(c+c\cdot n^c)$.
Then $\prod_j\xi_j$ is again polynomially admissible.
\end{enumerate}\end{theorem}
We defer to separate work constructing and analyzing a generic representation $\dual{\xi}$ 
for the convex compact space space $\dual{X}=\LIP(X,[0;1])$ of 
non-expansive functions $f:X\to[0;1]$, equipped with the supremum metric.

\subsection{Deferred Proofs from Section~\ref{s:Admissible2}}
\label{ss:Proofs1}

The proof of Theorem~\ref{t:Standard} resembles that of Example~\ref{x:Dyadic}, with four generalizations: 
First, the sequence $\vec a_n\in\TWO^*$ of binary integer numerators to dyadic approximations up to error $2^{-n}$
is replaced by a sequence of indices (w.r.t. some arbitrary but fixed enumeration) of centers of closed balls of radius $2^{-n}$ covering $X$;
indices $\vec a_n\in\TWO^m$ whose binary lengths $m=|\vec a_n|$ correspond to the binary entropy $m=\Entropy(n)$ of $X$, by the definition of entropy.
(For each precision $n$, we only consider binary strings of same length $m$: to avoid dealing with delimiters
and since, for any range $M=M(n)\subseteq\{0,1,\ldots,m-1\}$ of lengths, $\TWO^M$ can be replaced by $\TWO^m$.)
This results in a representation with modulus of continuity $n\mapsto\sum_{m\leq n} \Entropy(m)$.
Secondly, approximation error $2^{-n-1}$ twice as good as the required $2^{-n}$ yields admissibility:
achieved in the real setting using `rounding to nearest' and convexity in Equation~\eqref{e:4Admissible4},
while in our abstract space $X$ requires explicitly proceeding from covering 
by closed balls of radius $2^{-n}$ to ones of radius $2^{-n-1}$, increasing the modulus of continuity to $n\mapsto\sum_{m\leq n} \Entropy(m+1)$.
Thirdly, inspired by the signed binary representation, 
instead of each subsequent approximation $x_{m}\in X$ for $m:=n+1$ entirely superceding the previous $x_n$ as in the dyadic representation,
consider a covering with radius $2^{-n-2}$ of $\cball(x_n,1/2^n)$ rather than of entire $X$,
which replaces the absolute entropy $\Entropy(m+1)$ by the possibly smaller relative $\Entropy_{m-1}(m+1)$:
yielding modulus of continuity $n\mapsto\sum_{m\leq n} \Entropy_{m-1}(m+1)$.
And finally, instead of merely incrementing the precision $n\mapsto n+1$ between subsequent approximations,
allow for `jumps' by considering any strictly increasing sequence $n_k$:
yielding modulus $n_K\mapsto\Entropy(1+n_0)+\sum_{k=0}^K \Entropy_{n_{k}}\big(1+n_{k+1}\big)$.

\begin{proof}[Theorem~\ref{t:Standard}]
We construct $\xi_K:\subseteq\TWO^{\mu(n_K)}\to X$ by induction on $K$ such that 
\begin{enumerate}
\item[I)] the balls $\cball_{\pmb{1+}n_{K+1}}\big(\xi_{K+1}(\vec a)\big)$, $\vec a\in\dom(\xi_{K+1})$, cover $\cball_{n_K}\big(\xi_K(\vec b)\big)$.
\item[II)] $\xi_{K+1}$ satisfies 
\[ \forall \vec a\in\dom(\xi_{K+1}): \quad 
 d\big(\xi_{K+1}(\vec a),\xi_{K}(\vec b)\big)\leq2^{-n_K}  \]  
\end{enumerate}
with the abbreviation $\vec b:=\vec a|_{<\mu(n_K)}\in\dom(\xi_K)$.

\medskip
For $K=0$, by definition of $\Entropy$, at most $2^{\Entropy(1+n_0)}$ closed balls of radius $2^{-n_0\pmb{-1}}$ suffice to cover $X$.
This cardinality bound asserts that their centers can be enumerated over $\TWO^{\Entropy(1+n_0)}$ by some $\xi_0$.
Note $\mu(n_0)=\Entropy(1+n_0)$.

For the induction step $K\mapsto K+1$, 
consider any $\vec a\in\dom(\xi_K)\subseteq\TWO^{\mu(n_K)}$.  
By definition of $\Entropy$, at most $2^{\Entropy_{n_K}(1+n_{K+1})}$ closed balls of radius $2^{-n_{K+1}\pmb{-1}}$ suffice 
to cover $\cball_{n_K}\big(\xi_K(\vec a)\big)$.
Let 
\[ \xi_{K+1}(\vec a\,\cdot):\;\subseteq\; \TWO^{\Entropy_{n_K}(1+n_{K+1})}\;\ni\;\vec b\;\mapsto\; 
   \xi_{K+1}(\vec a\,\vec b)\;\in\; \cball_{n_K}\big(\xi_K(\vec a)\big) \]
enumerate such a family of centers. Then (II) holds by construction; and
(I) holds by induction hypothesis combined with `transitivity' of coverings.
This concludes the induction step.

\medskip
From (I) and (II) it follows that the partial mapping, 
\[ \xi:\subseteq\Cantor\;\ni\;\bar a \;\mapsto\; \lim\nolimits_{K\to\infty} \xi_K\big(\bar a|_{<\mu(n_K)}\big) \;\in\; X \]
defined whenever the limit exists,
is surjective (I) and has modulus of continuity $\mu$ (II). 
The prerequisite asserts it to satisfy condition (b\,i)/(d\,i) from Definition~\ref{d:Admissible2}, respectively.

In order to establish condition (c\,i)/(e\,i) from Definition~\ref{d:Admissible2},
consider $x,x'\in X$ with $d(x,x')\leq2^{-n_K-1}$.
Since in (I) balls cover $X$ of half the radius required in (II),
there exists $\vec a\in\dom(\xi_K)$ with $d\big(\xi_K(\vec a),x\big)\leq2^{-n_K\pmb{-1}}$:
which the above construction extends to a $\xi$-name $\bar a$ of $x$;
as well as to a $\xi$-name $\bar a'$ of $x'$, due to $d\big(\xi_K(\vec a),x'\big)\leq2^{-n_K}$.
Hence $\bar a$ and $\bar a'$ agree on the first $\mu(n_K)$ symbols:
$d_\xi(x,x')\leq2^{-\mu(n_K)}\leq 2^{-\Entropy(n_K)}$ by Lemma~\ref{l:Entropy}g).
In case $n\neq n_K$, proceed to the next larger $n_K\leq\calO(n)$ 
or to $n_K\leq\poly(n)$ according to the hypothesis.

In order to establish condition (b\,ii')/(d\,ii') from Definition~\ref{d:Admissible2},
fix some not necessarily surjective $\xi':\subseteq\Cantor\to X$ with minimal modulus of continuity $\mu'$,
and re-use the notation $\vec\Xi_m$ with
Properties~\eqref{e:4Admissible1},\eqref{e:4Admissible2},\eqref{e:4Admissible3}.
As replacement for Equation~\eqref{e:4Admissible4}, let $r(\vec u)$ denote any element of $\xi'[\vec u\circ\Cantor]$.
Now use (I) and (II) to construct, similarly to $\xi=\lim_K \xi_K$ by induction on $K\in\IN$, 
some mapping $F_K:\vec\Xi_{\mu'(n_K\pmb{+1})}\to\dom(\xi_K)$ such that
\[ \forall \vec u \in\vec\Xi_{\mu'(n_K\pmb{+1})}: \quad
r(\vec u)\;\in\;\cball\Big(\xi_K\big(F_K(\vec u)\big),2^{-n_K\pmb{-1}}\Big) \enspace . \]
Then, by triangle inequality and the definition of $\mu'$,
\begin{equation}
\label{e:4Admissible7}
\xi'[\vec u\circ\Cantor]\;\subseteq\;\cball\big(r(\vec u),2^{-n_K\pmb{-1}}\big)
\;\subseteq\;\cball\Big(\xi_K\big(F_K(\vec u)\big),2^{-n_K}\Big) \enspace : 
\end{equation}
Hence 
$ F:\dom(\xi')\;\ni\;\bar u\;\mapsto\;\lim\nolimits_{K\to\infty} F_K\big(\bar u|_{<\mu'(n_K+1)}\big) \;\in\; \dom(\xi) $
is well-defined, maps any $\xi'$-name of any $x\in X$ to a $\xi$-name of the same $x$,
and has modulus of continuity $\nu:\mu(n_K)\mapsto\mu'(n_K+1)$.
The latter means $\nu\circ\mu(n)\leq\mu'(n+1)$ for $n=n_K$,
and for $n\neq n_K$ implies $\nu\circ\mu(n)\leq\mu'\big(\calO(n)\big)$ 
by proceeding to the next larger $n_K\leq\calO(n)$ in the linear case,
or $\nu\circ\mu(n)\leq\mu'\big(\poly(n)\big)$ by proceeding to $n_K\leq\poly(n)$ in the polynomial case.
\qed\end{proof}

\begin{proof}[Theorem~\ref{t:Main1}]
\begin{enumerate}
\item[a\,i)]
Since $\xi$ is polynomially standard,
any not necessarily surjective $\xi'':\subseteq\Cantor\to X$ has
$\xi''\reduceqP\xi$ by Definition~\ref{d:Admissible2}(b\,ii'); 
and since $\xi'$ is polynomially admissible and $\xi$ is surjective,
Definition~\ref{d:Admissible2}(b\,ii) implies $\xi\reduceqP\xi'$.
Hence $\xi''\reduceqP\xi'$ by transitivity.
\\
To see that $\xi'$ satisfies Definition~\ref{d:Admissible2}(c\,i), 
exploit $\poly\circ\poly=\poly$ and
suppose $d(x,x')\leq2^{-\poly\circ\poly(n)}=2^{-\poly(m)}$ for $m:=\poly(n)$.
Since $\xi$ does satisfy Definition~\ref{d:Admissible2}(c\,i),
there exist $\xi$-names $\bar u$ of $x$ and $\bar u'$ of $x'$
with $D(\bar u,\bar u')=d_\xi(x,x')\leq2^{-\Entropy(m)}=2^{-\Entropy\circ\poly(n)}$.
Now $\xi\reduceqP\xi'$ yields $H:\dom(\xi)\to\dom(\xi')$
with modulus $\lambda$ such that $\xi\sqsubseteq\xi'\circ H$
and $\lambda\circ\mu'\leq\mu\circ\poly\leq\Entropy\circ\poly$
by Definition~\ref{d:Admissible2}(b\,i).
Finally note that, abbreviating $n':=\mu'(n)$,
$\lambda(n')\leq\Entropy\circ\poly(n)$
means that $H$ maps $\bar u$ and $\bar u'$ 
with $D(\bar u,\bar u')\leq2^{-\Entropy\circ\poly(n)}$
to $\bar v:=H(\bar u)$ and $\bar v':=H(\bar v')$
with $D(\bar v,\bar v')\leq 2^{-n'}=2^{-\mu'(n)}\leq 2^{-\Entropy(n)}$.
\item[a\,ii)] 
similarly.
\item[b\,i)]
Let $d(x,x')\leq2^{-\poly\circ\lambda(n)}=2^{-\poly(m)}$ for $m:=\lambda(n)$.
Since $\xi$ satisfies Definition~\ref{d:Admissible2}(c\,i),
there exist $\xi$-names $\bar u$ of $x$ and $\bar u'$ of $x'$
with $D(\bar u,\bar u')\leq2^{-\Entropy(m)}=2^{-\Entropy\circ\lambda(n)}$.
By hypothesis, $f\circ\xi$ maps these to $y:=f\circ\xi(\bar u)=f(x)$
and $y':=f\circ\xi(\bar u')=f(x')$ with distance $e(y,y')\leq 2^{-n}$.
\item[b\,ii)] 
similarly.
\item[c\,i)]
Let $\nu$ denote the minimal modulus of continuity of $\upsilon$.
Consider $\upsilon':=f\circ\xi:\subseteq\Cantor\to Y$
with modulus of continuity $\nu'\leq\mu\circ\lambda$. 
Definition~\ref{d:Admissible2}(b\,ii') implies $\upsilon'\reduceqP\upsilon$:
$f\circ\xi=\upsilon'=\upsilon\circ F$ for some $F:\dom(\upsilon')\to\dom(\upsilon)$
with modulus $\Lambda$ satisfying $\upsilon'=\upsilon\circ F$ and 
$\Lambda\circ\nu\leq\nu'\circ\poly$.
Now $\theta\leq\nu$ by Lemma~\ref{l:Entropy}a+g);
and $\mu\leq\Entropy\circ\poly$ by Definition~\ref{d:Admissible2}(b\,i).
\item[c\,ii)]
similarly.
\item[d\,i)]s
Consider $f\circ\xi\sqsubseteq\upsilon\circ F$
with modulus of continuity $\Lambda\circ\nu\leq(\Entropy\circ\lambda)\circ(\theta\circ\poly)$.
Now apply (b\,i).
\item[d\,ii)]
similarly.
\qed\end{enumerate}\end{proof}

\begin{proof}[Theorem~\ref{t:Cartesian}] 
\begin{enumerate}
\item[a)]
For $\bar u,\bar v\in\Cantor$ and non-decreasing $\mu,\nu:\IN\to\IN$, abbreviate
\begin{multline*}
\bar u  \myoplus{\mu}{\nu} \bar v \;:=\;
\big( u_0, u_1, \ldots, u_{\mu(0)-1} , \: v_0, v_1, \ldots, v_{\nu(0)-1} , \\
u_{\mu(0)}, u_{\mu(0)+1}, \ldots, u_{\mu(1)-1} , \: v_{\nu(0)}, v_{\nu(0)+1},  \ldots, v_{\nu(1)-1} , \:\: \ldots\ldots \\
u_{\mu(n)}, u_{\mu(n)+1}, \ldots, u_{\mu(n+1)-1} , \: v_{\nu(n)}, v_{\nu(n)+1}, \ldots, v_{\nu(n+1)-1} , \:\: \ldots \big) 
\end{multline*}
and similarly let $F\myoplus{\mu}{\nu} G$ pointwise for $F,G:\subseteq\Cantor\to\Cantor$.

Now let $\Entropy$ and $\theta$ denote the binary entropies of $(X,d)$ and $(Y,e)$, respectively;
and let $\mu_\xi,\mu_\upsilon$ denote the minimal moduli of continuity of $\xi$ and $\upsilon$, respectively. Then 
\[ \xi\times\upsilon \::\: (\bar u  \myoplus{\mu_\xi}{\mu_\upsilon} \bar v) \;\mapsto\; \big(\xi(\bar u),\upsilon(\bar v)\big) \;\in\; X\times Y \]
has modulus of continuity $\mu_\xi+\mu_\upsilon$.
We also record that, if $F$ and $G$ have moduli of continuity $\lambda_F$ and $\lambda_G$, respectively, then 
\[ H\;:=\;F \myoplus{\mu}{\nu} G \quad\text{has modulus}\quad \mu(n)+\nu(n)\mapsto \max\big\{\lambda_F\circ\mu(n),\lambda_G\circ\nu(n)\big\} \enspace . \]
\item[a\,i)]
By hypothesis, it holds $\mu_\xi\leq\Entropy\circ\calO$ and $\mu_\upsilon\leq\theta\circ\calO$.
Hence $\mu_\xi+\mu_\upsilon\leq (\Entropy+\theta)\circ\calO$ with 
$\Entropy(n+1)+\theta(n+1)+1$ 
an upper bound to the entropy of $X\times Y$ according to Lemma~\ref{l:Entropy}c).
\\
Let $\delta:\subseteq\Cantor\twoheadrightarrow X\times Y$ denote a representation with modulus $\lambda$.
Then $\Pi_1\circ\delta$ is a representation of $X$ with modulus $\lambda$,
where $\Pi_1:X\times Y\ni (x,y)\mapsto x\in X$ denotes the 1-Lipschitz projection on the first coordinate.
Hence $\Pi_1\circ\delta\reduceqL\xi$, since $\xi$ is linearly admissible:
$\xi\circ F=\Pi_1\circ\delta$ for some $F:\dom(\delta)\to\dom(\xi)$ 
whose modulus $\lambda_F$ satisfies $\lambda_F\circ\mu_\xi\leq\lambda\circ\calO$.
Similarly $\upsilon\circ G=\Pi_2\circ\delta$ for some $G:\dom(\delta)\to\dom(\upsilon)$ 
for some $G:\dom(\delta)\to\dom(\upsilon)$ 
whose modulus $\lambda_G$ satisfies $\lambda_G\circ\mu_\upsilon\leq\lambda\circ\calO$. Then it holds
\[ (\xi\times\upsilon)\;\circ\; H\;=\;\delta \quad\text{for}\quad H\;:=\;F \myoplus{\mu_\xi}{\mu_\upsilon} G\]
and $H$ has modulus 
\[ \mu_\xi(n)+\mu_\upsilon(n) \;\mapsto\;
\max\big\{\lambda_F\circ\mu_\xi(n),\lambda_G\circ\mu_\upsilon(n)\big\}  \;\leq\; \lambda\circ\calO(n) \]
\item[a\,ii)]
Let $x,x'\in X$ with $d(x,x')\leq2^{-\calO(n)}$;
and let $y,y'\in Y$ with $e(y,y')\leq2^{-\calO(n)}$.
By hypothesis there exist $\bar u,\bar u',\bar v,\bar v'$
with $\xi(\bar u)=x$, $\xi(\bar u')=x'$, $\upsilon(\bar v)=y$, $\upsilon(\bar v')=y'$
such that $D(\bar u,\bar u')\leq 2^{-\Entropy(n)}$ and $D(\bar v,\bar v')\leq 2^{-\theta(n)}$.
Then $\bar w:=\bar u \myoplus{\mu_\xi}{\mu_\upsilon} \bar v$ and $\bar w':=\bar u' \myoplus{\mu_\xi}{\mu_\upsilon} \bar v'$
have $D(\bar w,\bar w')\leq 2^{-(\Entropy(n)+\theta(n))}$.
And replacing $n$ with $n'=\calO(n)$ bounds $\Entropy(n')+\theta(n')$ by the entropy of $X\times Y$ 
according to Lemma~\ref{l:Entropy}c).
\item[a\,iii)] and (a\,iv) similarly.
\item[b)]
Let $\mu_j:\IN\to\IN$ be non-decreasing and unbounded for at least one $j\in\IN$.
Abbreviate $\bar\mu(n):=\mu_0(n) + \mu_1(n-1) + \cdots + \mu_n(0)$ and observe
\begin{multline*}
\bar\mu(n+1)-\bar\mu(n) \;=\;
\big(\mu_0(n+1)-\mu_0(n)\big) \:+\:
\big(\mu_1(n)-\mu_1(n-1)\big) \:+\: \cdots \:+\: \\
+ \:
\big(\mu_j(n-j+1)-\mu_j(n-j)\big) \:+\: \cdots \:+\:
\big(\mu_n(1)-\mu_n(0)\big) \:+\: \mu_{n+1}(0)  \enspace .
\end{multline*}
Therefore there exists a bijective map $\myOplus{\mu_j}{j\in\IN}=\bar\imath:\Cantor^\omega \rightarrow \Cantor$ 
such that for $\bar v = \myOplus{\mu_j}{j\in\IN} \bar u_j=\bar\imath \big((\bar u_j)_{_j}\big)$,
the finite substring $\big(v_{\bar\mu(n)},\ldots,v_{\bar\mu(n+1)-1}\big)$
is the concatenation of finite (possibly empty) substrings
\begin{multline*}
\big(u_{0,\mu_0(n)},\ldots,u_{0,\mu_0(n+1)-1}\big), \quad
\big(u_{1,\mu_1(n-1)},\ldots,u_{1,\mu_1(n)-1}\big), \quad\ldots\\
\big(u_{j,\mu_j(n-j)},\ldots,u_{1,\mu_j(n-j+1)-1}\big), \quad\ldots\\
\big(u_{n,\mu_n(0)},\ldots,u_{n,\mu_n(1)-1}\big), \quad
\big(u_{n+1,0},\ldots,u_{n+1,\mu_{n+1}(0)-1}\big) \enspace :
\end{multline*}
Note that the $j$-th component of the inverse
$\bar\imath^{-1}_j:\Cantor \rightarrow \Cantor$ 
has modulus of continuity
$n\mapsto\mu\circ\loinv{\mu_j}(n+j)$, that is,
$\mu_j(n-j)\mapsto\bar\mu(n)$.
Extend $\myOplus{\mu_j}{j\in\IN} F_j=\bar\imath (F_j)_{_j}$ pointwise to $F_j:\subseteq\Cantor\to\Cantor$.

Now take as $\mu_j$ the minimal moduli of continuity of $\xi_j:\subseteq\Cantor\twoheadrightarrow X_j$
and consider the representation 
\begin{multline*}
 \bar\xi\;=\;\prod\nolimits_j \xi_j \;:=\; \prod\nolimits_j (\xi_j\circ\bar\imath^{-1}_j) \;\subseteq\; \Cantor \\
\;\ni\; \myOplus{\mu_j}{j\in\IN} \bar u_j 
\;\mapsto\; \big(\xi_j(\bar u_j)\big)_j \;\in\; \prod\nolimits_j X_j 
\end{multline*}
having modulus of continuity $\bar\mu(n):=\sum\nolimits_{j \le n} \mu_j(n-j)$.
We also record that, if $F_j$ has modulus of continuity $\lambda_j$, then
\[ \myOplus{j\in\IN}{\mu_j} F_j \quad\text{has modulus}\quad
\sum\nolimits_{j \le n} \mu_j(n-j)\;\mapsto\; \max\nolimits_j \lambda_j\circ\mu_j(n) \enspace . \]
Let $\Entropy_j$ denote the entropy of $X_j$ and $\bar\Entropy$ the entropy of $\prod_j X_j/2^j$,
satisfying $\bar\Entropy(n)+1 \leq\sum\nolimits_{j<n} \Entropy_j\big(C+C\cdot(n-j+1)\big)$
by Lemma~\ref{l:Entropy}d) for some sufficiently large $C>0$ since $\diam(X_j)\geq1/2^j$ implies $\Entropy_j(j+2)\geq1$.
\item[b\,i)]
It follows $\bar\mu(n)=$
\[ \sum\nolimits_{j \le n} \mu_j(n-j) 
\;\le\; \sum\nolimits_{j \le n} \Entropy_j\big(c + c \cdot (n-j)\big) 
\;\leq\; \bar\Entropy\big(c + c \cdot (n+1)\big) \:+1 
\;\leq\; \bar\Entropy\big(\calO(n)\big) \] 
by hypothesis of $\xi_j$ being uniformly linearly admissible.
\\
Let $\delta:\subseteq\Cantor\twoheadrightarrow \prod_j X_j$ denote a representation with modulus $\lambda$.
Then $\Pi_j\circ\delta$ is a representation of $X_j$ with modulus $\lambda$.
Hence, by hypothesis of $\xi_j$ being uniformly linearly admissible, it holds $\Pi_j\circ\delta\reduceqL\xi_j$:
$\xi_j\circ F_j=\Pi_j\circ\delta$ for some $F_j:\dom(\delta)\to\dom(\xi_j)$ 
whose modulus $\lambda_j$ satisfies $\lambda_j\circ\mu_j(n)\leq\lambda(c+c\cdot n)$. Then 
\[ \big(\prod\nolimits_j \xi_j\big)\circ\big(\myOplus{j\in\IN}{\mu_j} F_j\big) \;=\;\delta \]
and $\myOplus{j\in\IN}{\mu_j} F_j$ has modulus 
\[ \bar\mu(n) \;=\; \sum\nolimits_{j \le n} \mu_j(n-j)\;\mapsto\; \max\nolimits_j \lambda_j\circ\mu_j(n) \;\leq\; \lambda\circ\calO(n) \]
\item[b\,ii)]
Let $(x_j)_{_j},(x_j')_{_j}\in \prod_j X_j$ with $d_j(x_j,x_j')\leq2^{-c-c(n-j)}$.
By hypothesis there exist $\bar u_j,\bar u_j'\in\Cantor$
with $\xi_j(\bar u_j)=x_j$ and $\xi_j(\bar u_j')=x_j'$
such that $D(\bar u_j,\bar u_j')\leq 2^{-\Entropy_j(n-j)}$.
Then $\bar w:=\myOplus{j\in\IN}{\mu_j} \bar u_j$ and $\bar w':=\myOplus{j\in\IN}{\mu_j} \bar u'_j$
have $D(\bar w,\bar w')\leq 2^{-\sum_{j<n}\Entropy_j(n-j)}$.
And replacing $n$ with $n'=\calO(n)$ bounds $\sum_{j<n'}\Entropy_j(n'-j)$ by $\bar\Entropy(n)$
according to Lemma~\ref{l:Entropy}c).
\item[b\,iii)] and (b\,iv) similarly. 
\qed\end{enumerate}\end{proof}

\section{Admissibility via Continuity for Multifunctions}
\label{s:Multifunc}

An admissible representation $\xi$ is usually not injective, 
hence $\xi^{-1}$ constitutes a \emph{multi}function.
Qualitative admissibility of continuous $\xi$ (Definition~\ref{d:Admissible})
could now be rephrased as follows: For every continuous $\xi'$, 
the \emph{multi}function $\xi^{-1}\circ\xi'$ admits a continuous selection $G$;
a H\"{o}lder-continuous one for linearly standard (Definition~\ref{d:Admissible2}).
Similarly, continuous $f:X\to Y$ having a continuous $(\xi,\upsilon)$-realizer
(Fact~\ref{f:KreitzWeihrauch}) could be rephrased as follows:
the \emph{multi}function $\upsilon^{-1}\circ f\circ\xi$
admits a continuous selection $F$;
a H\"{o}lder-continuous one in case of Corollary~\ref{c:Main1}ii).

In this section we adapt from \cite{PZ13} a logical notion of continuity for multifunctions
that (i) generalizes from the single-valued case, (ii) is closed under composition and restriction,
and (iii) yields a quantitatively continuous selection provided that it maps from/to \emph{ultra}metric spaces.
Admissibility of $\xi$ thus elegantly means that both $\xi$ and $\xi^{-1}$ are continuous;
and also the \emph{Main Theorem} becomes a mere consequence of (i)+(ii)+(iii).

\subsection{Multifunctions, Restriction and Selection}
\label{ss:Multifunc}

A multifunction is an oxymoron, namely a function with\emph{out} the axiom of functionality aka extensionality:
\[ x=x', \; y=f(x), \; y'=f(x') \quad \Rightarrow \quad y=y' \enspace .\]
Such generalization is unavoidable in real computation \cite{Luc77}.
For example, the Heaviside function $h_0$ is discontinuous and hence uncomputable,
yet for $\varepsilon>0$ its `soft'/non-extensional/multi-functional variant is computable:
\begin{equation}
\label{e:Heaviside}
h_\varepsilon(t) \;:=\; \left\{ \begin{array}{rl} 0 &: t\leq \varepsilon \\ 1&: t\ge-\varepsilon 
\end{array} \right. 
\end{equation}
Formally, a partial multivalued function (multifunction)
$F$ between sets $X,Y$ is a relation $F\subseteq X\times Y$
that models a computational \emph{search} problem:
Given (any name of) $x\in X$, return some (name of some) $y\in Y$ with $(x,y)\in F$.
Mathematically one may identify the relation $f$ with the single-valued total function
$F:X\ni x\mapsto \{y\in Y \mid (x,y)\in F\}$
from $X$ to the powerset $2^Y$; but we prefer 
the notation $f\pcolon X\toto Y$ to emphasize
that not every $y\in F(x)$ needs to occur as output;
see also Example~\ref{x:Continuity2}.
Letting the answer $y$ depend on the code of $x$
means dropping the requirement for ordinary functions to be extensional;
hence, in spite of the oxymoron, such $F$ is also called a \emph{non-}extensional function.
Note that no output is feasible in case $F(x)=\emptyset$.

\begin{definition}
\label{d:Multifunc}
Abbreviate with $\dom(F) \dfeq \{x \mid F(x)\neq\emptyset\}$
the domain of $F$; and $\range(F)\dfeq \{y\mid \exists x: (x,y)\in F\}$.
$F$ is \emph{total} in case $\dom(F)=X$;
\emph{surjective} in case $\range(F)=Y$.
The composition of partial multifunctions $F:\subseteq X\toto Y$ 
and $G:\subseteq Y\toto Z$ is $G\circ F=$
\[ 
  \big\{(x,z) \;\big|\; x\in X, z\in Z, F(x)\subseteq\dom(G),
  \;\exists y\in Y: (x,y)\in F\wedge (y,z)\in G \} 
\] 
$F$ is \emph{compact} if the image $F[C]\subseteq Y$ is compact for every compact $C\subseteq\dom(F)$.
\end{definition}
Note that Definition~\ref{d:Multifunc} indeed generalizes from the single-valued case.

\begin{proposition}
\label{p:Continuity}
\begin{enumerate}
\item[a)]
Suppose $(X,d)$ and $(Y,e)$ are compact and $f:X\to Y$ is single-valued.
Then $f$ is continuous iff it is (i) compact and (ii) $f^{-1}(y)\subseteq X$ is closed for every $y\in Y$.
\item[b)]
For a continuous total single-valued function $f:(X,d)\to (Y,e)$ with compact domain,
its multivalued partial inverse $f^{-1}:\subseteq Y\toto X$ is compact.
\item[c)]
Suppose $(X,d)$ and $(Y,e)$ are compact.
A partial multifunction $F:\subseteq X\toto Y$ and its partial inverse $F^{-1}:\subseteq Y\toto X$
are both simultaneously compact iff the graph $F\subseteq X\times Y$ is a compact set.
\item[d)]
If both total multifunctions $F:X\toto Y$ and $G:Y\toto Z$ are compact,
then so is their composition $G\circ F$.
\end{enumerate}
\end{proposition}
\begin{proof} 
\begin{enumerate}
\item[a)] 
$(\Rightarrow)$ Omitted.
$(\Leftarrow)$ Assume conditions (i) and (ii).
We may also assume that $f$ is surjective.
Fix $a \in X$ and the family $\mathcal{H}$ of 
closed neighborhoods of $a$.
For every $x \neq a$, 
there exists $H \in \mathcal{H}$ disjoint to $f^{-1}(f(x))$.
This shows $\bigcap_{H \in \mathcal{H}} f[H] = \{f(a)\}$.
Suppose that $f$ is discontinuous at $a$.
There exists an open neighborhood $U$ of $f(a)$ such that
$f[H] \nsubseteq U$ for every $H \in \mathcal{H}$.
Then $\{U^c\} \cup \{f[H]\}_{H \in \mathcal{H}}$
form a collection of closed subsets in $Y$
satisfying finite intersection property but
$U^c \cap \bigcap_{H \in \mathcal{H}} f[H] = \emptyset$, 
a contradiction\footnote{See \textsc{Brian M. Scott} answering
\url{https://math.stackexchange.com/questions/1527612/}}.
\item[b)]
Observe that, in metric spaces, compactness implies closedness
and that in compact metric space compactness and closedness coincide.
\item[c)]
$(\Rightarrow)$ Consider an open covering $\mathcal{U}$ of
$F \subseteq X \times Y$. 
For each $x \in X$, we have a finite subcovering
$\mathcal{U}_x$ of $F \cap (\{x\} \times Y)$. 
Then there exists open $V_x \in X$ such that
$x \in V$ and 
$F \cap (V_x \times Y) \subseteq \bigcup \mathcal{U}_x$.
We can find a covering $\{ V_{x_1}, \cdots, V_{x_n} \}$ of $X$.
Then $\mathcal{U}_{x_1} \cup \cdots \cup \mathcal{U}_{x_n}$
is the desired finite subcovering.
$(\Leftarrow)$
For compact $C \subseteq X$, we have
$F[C] = \pi_Y[(C \times Y) \cap F]$ where
$\pi_Y:X\times Y \rightarrow Y$ is the projection.
\item[d)] 
We have $G \circ F(x) = G[F(x)]$
for $x\in X$ since $G$ is total.
\qed\end{enumerate}
\end{proof}
A computational problem, considered as total single-valued function $f:X\to Y$, 
becomes `easier' when \emph{restricting} to arguments $x\in X'\subseteq X$,
that is, when proceeding to $f'=f|_{X'}$ for some $X'\subseteq X$.
A search problem, considered as total multifunction $F:X\toto Y$, 
additionally becomes `easier' when proceeding to any $F'\subseteq X\toto Y$
satisfying the following: $F'(x)\supseteq F(x)$ for every $x\in\dom(F')$.
We call such $F'$ also a \emph{restriction} of $F$, and write $F'\sqsubseteq F$.
A single-valued function $f:\dom(F)\to Y$ is a \emph{selection} of $F:\subseteq X\toto Y$
if $F$ is a restriction of $f$.

\begin{observation}
\label{o:Restriction}
Fix partial multifunctions 
$F:\subseteq X\toto Y$ and $G:\subseteq Y\toto Z$.
\begin{enumerate}
\item[a)]
The composition of restrictions $F'\sqsubseteq F$ and $G'\sqsubseteq G$,
is again a restriction $G'\circ F'\sqsubseteq G\circ F$.
\item[b)]
It holds $F^{-1}\circ F\sqsubseteq\id_X:X\to X$.
Single-valued surjective partial $g:\subseteq X\twoheadrightarrow Y$ 
furthermore satisfy $g\circ g^{-1}=\id_Y$.
\item[c)]
Consider multifunction $f:X\toto Y$, representations $\xi:\subseteq U\twoheadrightarrow X$ 
and $\upsilon:\subseteq V\twoheadrightarrow Y$, and single-valued function $F:\subseteq U\to V$.
The following are equivalent: \\ 
(i) $f\circ\xi\sqsubseteq\upsilon\circ F$ \\
(ii) $f\sqsubseteq\upsilon\circ F\circ\xi^{-1}$ \\
(iii) $\upsilon^{-1}\circ f\circ\xi\sqsubseteq F$.
\end{enumerate}
\end{observation}

\subsection{Quantitative Continuity for Multifunctions}
\label{ss:Continuity}

(Lower/upper) \emph{hemi}continuity is a well-established generalization to multi-functions;
however, it violates closure under the above generalized conception of restriction.
Considering a multivalued $f$ as function with values in the hyperspace of compact subsets 
gives rise to a notion of continuity incompatible with computability; 
see Example~\ref{x:Continuity2} and Remark~\ref{r:Selection} below.
Instead Definition~\ref{d:Multicont} below adapts, and quantitatively refines, a notion of `sequential' uniform continuity 
for multifunctions \cite{BrHe94} \cite[\S4+\S6]{PZ13} \cite[Theorem~3.8]{Schroeder2020}---such as to satisfy the following properties and examples:

\begin{lemma}
\label{l:Continuity}
\begin{enumerate}
\item[a)]
A single-valued function has modulus of sequential continuity $\mu\big(n+\calO(1)\big)$ iff it is
$\mu\big(n+\calO(1)\big)$-continuous when considered as a multifunction. 
\item[b)]
Suppose that $F:\subseteq X\toto Y$ is $\mu$-continuous.
Then every restriction $F'\sqsubseteq F$ is again $\mu$-continuous.
\item[c)]
If additionally $G:\subseteq Y\toto Z$ is $\nu$-continuous,
then $G\circ F$ is $\mu\circ\nu$-continuous
\end{enumerate}
\end{lemma}
To motivate Definition~\ref{d:Multicont}, 
consider the following `sequential' weakening of Equation~\eqref{e:Modulus},
qualitatively still equivalent to uniform continuity:

\begin{observation}
\label{o:Continuity}
Fix $f:X\to Y$ and \emph{strictly} increasing $\mu:\IN\to\IN$ and 
some `offset' $n_0\in\IN$.
Equation~(\ref{e:Modulus}) holds for all $n\geq n_0$ 
iff the following is true for some/every positive $k\in\IN$:
\begin{multline} \label{e:Modulus3}
\forall x_0,\ldots,x_k\in X \; 
\forall n_1,\ldots,n_k\in\IN \\\
\Big( n_1\geq n_0 \;\wedge\; 
\bigwedge\nolimits_{j=2}^{k} n_{j}\geq n_{j-1}+n_0 
\;\wedge\; \bigwedge\nolimits_{j=1}^k d(x_{j},x_{j-1})\leq 2^{-\mu(n_{j})}  \\
\Longrightarrow\;\; \bigwedge\nolimits_{j=1}^k 
e(y_{j},y_{j-1})\leq 2^{-n_{j}}\Big)  \enspace, \qquad y_j=f(x_j) 
\end{multline}
\end{observation}
In this case let us call $\mu$ a modulus of \emph{sequential} continuity of $f$.
For $\nu$ a modulus of sequential continuity of $g:Y\to Z$
(with offset $m_0$),
$\mu\circ\nu$ is a modulus of sequential continuity of $g\circ f$:
with offset $\max(n_0,m_0)$ instead of $n_0$,
due to strict monotonicity of $\nu$.
Every \Martin{strictly increasing} modulus of continuity of a function is also one of sequential continuity;
and modulus $\mu=\mu(n)$ of sequential continuity conversely 
yields a modulus $\mu\big(n+\calO(1)\big)$ of continuity
in the original sense of Equation~\eqref{e:Modulus}.
We have thus motivated the following definition:

\begin{definition}
\label{d:Multicont}
Fix metric spaces $(X,d)$ and $(Y,e)$ and strictly increasing $\mu:\IN\to\IN$.
A total multifunction $F:X\toto Y$
is called \emph{$\mu$-continuous} if there exists some $n_0\in\IN$,
and to every $x_0\in X$ there exists some $y_0\in F(x_0)$,
such that the following holds for every positive $k\in\IN$:
\begin{multline}
\label{e:Multicont}
\forall n_1\geq n_0 \;\; \forall x_1\in\cball_{\mu(n_1)}(x_0) \;\; \exists y_1\in F(x_1)\cap\cball_{n_1}(y_0)  \\  
\forall n_2\geq n_1+n_0 \;\; \forall x_2\in\cball_{\mu(n_2)}(x_1) \;\; \exists y_2\in F(x_2)\cap\cball_{n_2}(y_1) \qquad \ldots \quad \\  
\forall n_{k}\geq n_{k-1}+n_0 \;\; \forall x_{k}\in\cball_{\mu(n_{k})}(x_{k-1}) \;\;\exists y_{k}\in F(x_{k})\cap\cball_{n_{k}}(y_{k-1})
\end{multline}
with the abbreviation $\cball_n(x):=\cball(x,2^{-n})$.
\end{definition}
This notion is not a topological one \cite{ksc2020park}; and it exceeds Classical Logic 
in that the number $\calO(k)$ of quantifiers varies: in order to capture finite sequences of arbitrary length.
Examples~\ref{x:Continuity2}c+d) below work with $n_0=1$;
however Definition~\ref{d:Multicont} should be \emph{qualitatively} invariant under metric rescaling: to which end we allow arbitrary $n_0$.
Closure under composition (Lemma~\ref{l:Continuity}c) relies on the modulus to be strictly increasing.
Parameter $n_0$ requires $n_1\geq n_0$ and $n_k\geq n_{k-1}+n_0\geq k n_0$ by induction, and hence for $k\geq\ell$:
\begin{gather}
\label{e:Multicont2}
d(x_\ell,x_k)\;\leq\; \sum\nolimits_{j=\ell+1}^k d(x_j,x_{j-1}) \;\leq\; 
\sum\nolimits_{j=\ell+1}^k 2^{-\mu(n_j)} \;\leq \\
\leq\; \sum\nolimits_{j=\ell+1}^k 2^{-n_j} \;\leq\;  \nonumber
\sum\nolimits_{j=\ell+1}^k 2^{-j n_0} \;\leq\; 2^{-n_0(\ell+1)}/(1-2^{-n_0}) \;\leq\; 2^{-(\ell+1)n_0+1} 
\end{gather}
as exploited in Example~\ref{x:Continuity2}; compare Observation~\ref{o:Continuity}.

\begin{example}
\label{x:Continuity2}
Recall the `soft' Heaviside function $h_\varepsilon$ from Equation~\eqref{e:Heaviside}.
\begin{enumerate}
\item[a)]
For every $\varepsilon>0$, the multifunction $h_\varepsilon:[-1;1]\toto[0;1]$ is $\id$-continuous;
and so is $\bar h:=1-h_\varepsilon:[-1;1]\toto[0;1]$;
but none is so for $\varepsilon=0$: in agreement with both being computable for $\varepsilon>0$
but not for $\varepsilon=0$.
\item[b)]
Consider $h_\varepsilon$ as (single-valued) function $h_\varepsilon:[0;1]\to\calK([0;1])$
with the hyperspace $\calK([0;1])$ from Lemma~\ref{l:Entropy}e).
For $\varepsilon>0$ it is dicontinuous at both $t=-\varepsilon$ and at $t=+\varepsilon$:
in disagreement with its in/computability.
\item[c)]
The multivalued inverse of the dyadic representation $\delta^{-1}$
is $(n+2)$-continuous with $n_0=1$.
\item[d)]
The multivalued inverse of the signed-digit representation $\sigma^{-1}$
is $(n+2)$-continuous with $n_0=1$.
\end{enumerate}
\end{example}
%

\subsection{Quantitative Selection on Ultrametric Spaces}
\label{ss:Selection}

Unlike continuous multifunctions over the reals \cite[Fig.5]{PZ13},
those over Cantor space do admit a continuous selection---and even a bound on its modulus:

\begin{theorem}
\label{t:Selection}
Suppose $(\calU,D)$ and $(\calV,E)$ are compact ultrametric spaces. 
If total $G:\calU\toto\calV$ is $\mu$-continuous and compact,
then $G$ admits a selection with modulus of continuity $\mu\big(n+\calO(1)\big)$;
with modulus $\mu$ in case $\calU$ is (a compact subset of) Cantor space $\Cantor$.
\end{theorem}
Recall Lemma~\ref{l:Ultrametric} for examples.

\begin{remark}
\label{r:Selection}
\begin{enumerate}
\item[a)]
The literature contains several notions of continuity for multifunctions \cite{BrHe94}.
We have mentioned \emph{hemi}continuity at the beginning of this subsection,
as well as its deficiency. 
Defining continuity of a multifunction via having a continuous realizer
makes the \emph{Main Theorem} a tautology, but constitutes an extrinsic
(rather than intrinsic) conception.
Example~\ref{x:Continuity2}b) exhibits
the difference between considering a multifunction $f:X\toto Y$ as computational search problem 
and as function with set-values: if the co-domain $\calK(Y)$ is equipped with the Hausdorff metric.
On the other hand equipping it with the \emph{lower Vietoris topology} agrees with computability \cite{OvertChoice}:
\item[b)]
The literature contains many selection results for 
the hyperspace $\calK(Y)$ of non-empty compact subsets of some metric space $(Y,e)$
\cite{Selection,Repovs2014}:
with respect to qualitative continuity---while our
Theorem~\ref{t:Selection} is quantitative.
\\
Recall that the lower Vietoris topology is induced by the quasi (\ie non-symmetric) metric 
$E(V,V')=\sup_{v\in V} \inf_{v'\in V'} e(v,v')$, the so-called `half' of the Hausdorff metric.
With respect to this $E$, 
a single-valued function $F:X\to\calK(Y)$ 
has modulus of continuity $\mu$ iff:
For some $n_0\in\IN$ and every $x_0\in X$ 
and \emph{every} $y_0\in F(x_0)$ 
and some/every $k\in\IN_+$,
Equation~(\ref{e:Multicont2}) holds.
Definition~\ref{d:Multicont} can thus be considered as replacing
the universal quantification over $y_0$ with an existential one:
in order to obtain closure under (both kinds of) restriction;
and our quantitative Theorem~\ref{t:Selection}
requires the underlying spaces to be ultrametric.
\item[c)]
The proof of Example~\ref{x:Continuity2}c) shows that $\delta^{-1}$
is actually $\calO(\sqrt{n})$-continuous---as one might hope 
since $\delta$ has minimal modulus $\calO(n^2)$ 
according to Example~\ref{x:Dyadic}iii).
However Definition~\ref{d:Multicont} can only consider strictly
increasing moduli, in order to ensure closure under composition:
This is a disadvantage of the present sections's approach to admissible representations via multivalued continuity.
On the other hand it yields a quantitative \emph{Main Theorem} 
for \emph{multi}functions $f:X\toto Y$, see Theorem~\ref{t:Main2}.
\end{enumerate}
\end{remark}

\subsection{Proof of the Selection Theorem}
\label{ss:Selection2}

In order to emphasize the in/dependencies among quantified variables,
let us rephrase Equation~\ref{e:Multicont} in Skolem form:

\begin{lemma}
\label{l:Continuity2}
Fix a strictly increasing $\mu:\IN\to\IN$ and 
a total compact multifunction $G:X\toto Y$ between metric spaces $(X,d)$ and $(Y,e)$.
$G$ is $\mu$-continuous iff there exists $n_0\in\IN$
and (single-valued but not necessarily continuous) partial functions $G_k:\subseteq X^{k+1}\times\IN^k\to Y$ 
satisfying the following for every $k\geq1$:
\begin{multline}
\label{e:Multicont3}
\forall x_0,x_1,\ldots,x_k\in X \quad
\forall n_1,n_2,\ldots,n_k\in\IN \\
\Big( n_1\geq n_0\wedge\bigwedge\nolimits_{j=2}^{k} n_j\geq n_{j-1}+n_0
\;\wedge\; \bigwedge\nolimits_{j=1}^k  d(x_j,x_{j-1})\leq 2^{-\mu(n_j)} \\
\Longrightarrow \quad
e\big(G_k(x_0,\ldots,x_k;n_1,\ldots,n_k),G_{k-1}(x_0,\ldots,x_{k-1};n_1,\ldots,n_{k-1})\big)\leq2^{-n_k} \Big)
\end{multline}
\end{lemma}
\begin{proof}
Equation~(\ref{e:Multicont3}) clearly implies (\ref{e:Multicont}).
For the converse, Equation~(\ref{e:Multicont}) yields 
partial Skolem functions $y_j=G_j^{(k)}(x_0,\ldots,x_j;n_1,\ldots,n_k)\in G(x_j)$, $0<j\leq k$,
satisfying Equation~(\ref{e:Multicont3}), but depending on $k$.
To remove said dependence on $k$, exploit compactness of $G$:
$y_j^{(k)}=G_j^{(k)}(x_0,\ldots,x_j;n_1,\ldots,n_k)\in G(x_j)$ admits a subsequence 
converging, as $k\to\infty$, to some $y_j=:G_j(x_0,\ldots,x_j;n_1,\ldots,n_k)\in G(x_j)$.
These thus pointwise, and by induction on $j$, defined functions $G_j$ satisfy Equation~(\ref{e:Multicont3}) 
for the infinitely many $k$ occurring in said subsequence---and therefore for all $k$.
\qed\end{proof}
\begin{figure}[htb]
\includegraphics[width=0.95\textwidth]{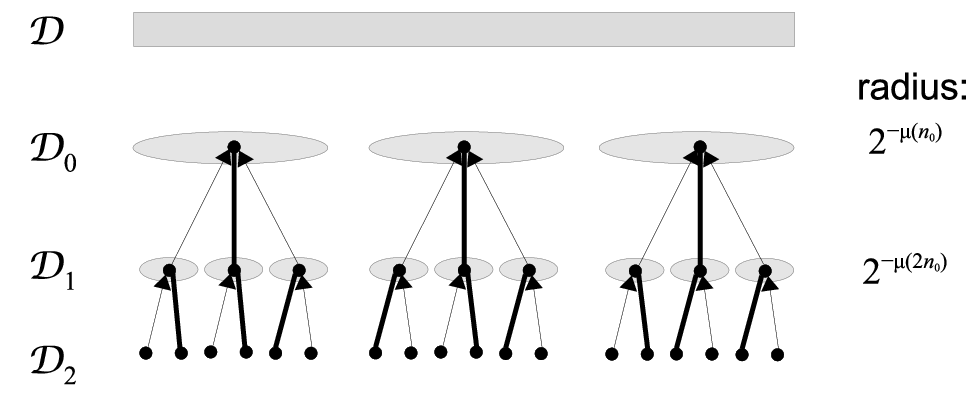}
\caption{\label{f:Ultradecomp}Hierarchical decomposition of the 
compact ultrametric space $(\calU,D)$ 
according to Lemma~\ref{l:Decomposition}.
(Fat dots connected by a bold line segment represent the same point in $\calU$.)}
\end{figure}

\begin{lemma}
\label{l:Decomposition}
Let $(\calU,D)$ denote a compact ultrametric space
and $\nu:\IN\to\IN$ be strictly increasing.
\\
Then $\calU$ admits a hierarchical cluster decomposition,
namely a non-decreasing sequence $\calU_m\subseteq\calU_{m+1}\subseteq\calU$ 
of finite sets with the following properties:
\begin{enumerate}
\item[a)]
The balls $\cball_{\nu(0)}(x)$, $x\in\calU_0$, 
are pairwise disjoint and cover $\calU$.
More precisely:
\begin{equation}
\label{e:Decomposition1}
\text{any two distinct $x,x'\in\calU_{0}$ ~have~ $D(x,x')>2^{-\nu(0)}$.}
\end{equation}
\item[b)]
For each $m\in\IN$ and $x\in\calU_m$,
the balls $\cball_{\nu(m+1)}(x')$, 
$x'\in\cball_{\nu(m)}(x)\cap\calU_{m+1}$,
partition $\cball_{\nu(m)}(x)$. More precisely, for any distinct $x',x''\in\calU$,
\begin{equation}
\label{e:Decomposition2}
x',x''\in\cball_{\nu(m)}(x)\cap\calU_{m+1} 
\quad\text{implies}\quad 2^{-\nu(m+1)}<D(x',x'')\leq 2^{-\nu(m)}
\end{equation}
\end{enumerate}
\end{lemma}
Such a decomposition can be constructed by inductive application of Observation~\ref{o:Ultra}:
pick and include in $\calU_{m+1}$ precisely one representative center $x'$ 
of all mutually equal balls $\cball_{\nu(m+1)}(x')$, $x'\in\cball_{\nu(m)}(x)$;
see Figure~\ref{f:Ultradecomp}.

\begin{proof}[Theorem~\ref{t:Selection}]
Abbreviate $\nu(m)=\mu(n_0+m\cdot n_0)$ for $m\in\IN$.
Consider the hierarchical cluster decomposition $\calU_m\subseteq\calU_{m+1}\subseteq\calU$ 
from Lemma~\ref{l:Decomposition}.
For $x\in\calU_{m+1}\setminus\calU_{m}$ let $\calP(x)\in\calU_m$ denote its \emph{predecessor},
that is, the unique $x'=\calP(x)\in\calU_m$ with 
$2^{-\nu(m+1)}<D(x,x')\leq 2^{\nu(m)}$.
This mapping $\calP$ extends to $\bigcup_m\calU_m\setminus\calU_0$.
Note that $\calP(\calP(x))$ may be undefined, namely in case $\calP(x)\in\calU_0$.
For $x\in\calU_{m+1}\setminus\calU_{m}$, let 
\begin{equation}
\label{e:Sequence}
\vec\calP(x)  \;:=\; \big(\,\calP^{(k)}(x)\,,\,\calP^{(k-1)}(x)\,,\,\ldots\,,\,\calP(\calP(x))\,,\,\calP(x),x\,\big) 
\end{equation}
abbreviate the sequence of iterated predecessors, reversed.
Here $|\vec\calP(x)|=k\leq m+1$ is the unique natural number such that $k$-fold application $\calP^{(k)}(x)\in\calU_0$.
Again note that the sequence (\ref{e:Sequence}) may `skip' some levels,
\ie, the length $k$ may be (much) shorter than $m$.
To $\vec x:=\vec\calP(x)$ choose $(n_1,\ldots,n_k)=\vec n(x)\in\IN$ maximal with $D(x_j,x_{j-1})\leq 2^{-\mu(n_j)}$:
Then Equation~(\ref{e:Decomposition1}) implies $n_1\geq n_0$ and $n_j\geq n_{j-1}+n_0$ for $j=2,\ldots k$;
hence the prerequisite of Equation~(\ref{e:Multicont3}) is satisfied.

\medskip
Next we construct a sequence of total single-valued functions $g_m:\calU_m\to\calV$ with the following properties:
\begin{enumerate}
\item[i)] $g_{m+1}|_{\calU_m}=g_m$
\item[ii)] $g_m$ has modulus of continuity $n\mapsto\mu(n+n_0-1)$
\item[iii)] $g_m\subseteq G$ as graphs, \ie $g_m$ is a selection of the restriction $G|_{\calU_m}$.
\end{enumerate}
These $g_m$ thus converge uniformly to some total function $g\subseteq G$ on $\calU$ with modulus of continuity $\mu(n+n_0-1)$.

\medskip
To obtain $(g_m)$, employ the (discontinuous) Skolem functions $G_k:\subseteq\calU^{k+1}\times\IN^k\to\calV$ from Lemma~\ref{l:Continuity2}:
Define $g_0(x)=G_0(x)$ for $x\in\calU_0$; and inductively, 
for $x\in\calU_{m+1}\setminus\calU_{m}$, define $g_{m+1}(x)=G_k\big(\vec\calP(x),\vec n(x)\big)$ with $k=|\vec\calP(x)|$.
As recorded above, $\vec\calP(x)$ and $\vec n(x)$ together satisfy the prerequisite of Equation~(\ref{e:Multicont3});
hence $g_{m+1}$ is well-defined and satisfies (i) and (iii) by inductive construction.

\medskip
Concerning (ii), observe that $g_0$ has sequential modulus of continuity $\mu$:
simply because the prerequisite of Equation~(\ref{e:Modulus}) with $n\geq n_0$
is not satisfied according to Equation~(\ref{e:Decomposition1}).

Regarding the induction step $g_m\mapsto g_{m+1}$, 
first consider $x\in\calU_{m+1}\setminus\calU_{m}$ and $x:=\calP(x')$:
$2^{-\nu(m+1)}<D(x,x')\leq2^{-\nu(m)}=2^{-\mu(n_0+mn_0)}$ as recorded above implies
$E(y,y')\leq2^{-(m+1)n_0}$ for $y=g_{m+1}(x)=G_k\big(\vec\calP(x),\vec n(x)\big)$
and $y'=g_{m+1}(x')=G_{k-1}\big(\vec\calP(x'),\vec n(x')\big)$, $k=|\vec\calP(x)|$.

Next consider $x',x''\in\calU_{m+1}\setminus\calU_m$:
If $2^{-\nu(m+1)}<D(x',x'')\leq2^{-\nu(m)}=2^{-\mu((m+1)n_0)}$, then $2^{-\nu(m+1)}<D(x',x),D(x'',x)\leq2^{-\nu(m)}$ 
for $x:=\calP(x')=\calP(x'')\in\calU_m$ according to Equation~(\ref{e:Decomposition2});
hence $E(y',y'')\leq\max\{E(y',y),E(y,y'')\}\leq2^{-(m+1)n_0}$ for $y=g_{m+1}(x)$, $y'=g_{m+1}(x')$, $y''=g_{m+1}(x'')$ as before.
Note that $D(x',x'')\leq2^{-\mu(n)}$ in case $(m+2)n_0>n>(m+1)n_0$ does not imply $D(x',x),d(x'',x)\leq2^{-\mu(n)}$;
hence here we only obtain modulus $n\mapsto\mu(n+n_0-1)$ instead of $\mu(n)$.
Note that this case cannot occur over (a subset of) Cantor space $\Cantor$ as domain, though.

Finally suppose $x,x'\in\calU_{m+1}\setminus\calU_m$ have $2^{-\nu(m)}<D(x,x')\leq2^{-\nu(m')}$.
Then $x'':=\calP(x)\in\calU_m$ and $x''':=\calP(x')\in\calU_m$ have 
$D(x'',x''')=D(x,x')$ by Equation~(\ref{e:Decomposition2}) 
and Observation~\ref{o:Ultra}iv); and 
$2^{-\nu(m+1)}<D(x,x''),D(x',x''')\leq2^{-\nu(m)}$ as recorded above.
Hence $y:=g_{m+1}(x)$, $y':=g_{m+1}(x')$, $y'':=g_{m+1}(x'')$, $y''':=g_{m+1}(x''')$
have $E(y,y''),E(y',y''')\leq2^{-(m+1)n_0}$ as before;
and $E(y'',y''')\leq2^{-(m'+1)n_0}$ by induction hypothesis:
Together $E(y,y')\leq \max\{E(y,y''),E(y'',y'''),E(y''',y')\}\leq2^{-(m'+1)n_0}$:
concluding the final case of the induction step to verify (ii).
\qed\end{proof}

\begin{proof}[Example~\ref{x:Continuity2}]
\begin{enumerate}
\item[a)]
Let $n_0\in\IN$ such that $2^{-n_0+1}<\varepsilon$.
If $x_0\leq 0$ take $y_0:=0\in h_\varepsilon(x_0)$, 
otherwise take $y_0:=1\in h_\varepsilon(x_0)$:
Equation~\eqref{e:Multicont2} shows $x_k<\varepsilon$
and hence $y_k:=0\in h_\varepsilon(x_k)$.
\item[c)]
Let us call an integer sequence $(a_n)_{_{n\geq1}}$ a $\tildedelta$-name of $x\in[0;1]$
iff $\big(\bin_n(a_n):n\geq1\big)\in\Cantor$ is a $\delta$-name of $x$.
Note that (only) here, the index starts with $1$ instead of $0$.
Since $\mu(n)=n+2$, Equation~\eqref{e:Multicont2} improves:
\begin{equation}
\label{e:Multicont4}
d(x_\ell,x_k) \;\leq\; \sum\nolimits_{j=\ell+1}^k 2^{-(n_j+2)} \;\leq\; 2^{-n_{\ell+1}-1} \enspace .
\end{equation}
\\
To $x_0\in[0;1]$ consider the integer sequence $\bar y_0=(y^{(0)}_n)_{_n}=\big(\lfloor 2^n\cdot x_0\rceil_n\big)$.
Since $|y^{(0)}_n-x_0|\leq2^{-n-1}$, and with $|x_0-x_k|\leq2^{-n_1-1}$ according to Equation~\eqref{e:Multicont4}, 
the $n_1$-element sequence $(y^{(0)}_1,\ldots,y^{(0)}_{n_1})$ extends to $\tildedelta$-names of any $x_k$.
More precisely, $\bar y^{(1)}=(y^{(1)}_n)_{_n}:=\big(\lfloor 2^n\cdot x_1\rceil_n\big)$ for $n>n_1$ and $y^{(1)}_n:=y^{(0)}_n$ for $n\leq n_1$
constitutes a $\tildedelta$-name of $x_1$: agreeing on the first $n_1$ entries;
and now with the property $|y^{(1)}_n-x_1|\leq2^{-n-1}$ for all $n>n_1$;
which, together with $|x_1-x_k|\leq2^{-n_2-1}$ makes $(y^{(0)}_1,y^{(0)}_{n_1},y^{(1)}_{n_1+1},y^{(1)}_{n_2})$ 
extend to $\tildedelta$-names of any $x_k$, $k>1$. And so forth:
yielding a sequence of $\tildedelta$-names $\bar y^{(k)}\in\IN^\IN$ of $x_k\in\cball_{n_k+2}(x_{k-1})$ 
which agree on the first $n_k$ entries.
These in turn yield $\delta$-names $y_k\in\delta^{-1}(x_k)\subseteq\Cantor$ which
agree on the first $\nu(k_k):=\sum_{j<n_k} j=n_k\cdot(n_k+1)/2\geq n_k$ bits: recall Example~\ref{x:Dyadic}.
In particular $y_{k}\in \cball_{n_{k}}(y_{k-1})$.
\item[d)]
We basically iterate the proof of Example~\ref{x:SignDigit}(iii):
Consider any $\vec b^{(0)}=\big(b_0^{(0)},\ldots,b_{n_1-1}^{(0)}\big)\in\THREE^{n_1}$ with $|x_0-\sbin(\vec b^{(0)})|\leq2^{-n_1-1}$:
such that it extends to signed binary expansions of any $x_k$,
since $|x_0-x_k|\leq2^{-n_1-1}$ according to Equation~\eqref{e:Multicont4}.
In particular the signed binary expansion of $\tfrac{1}{2}+(x_1-x_0)\cdot2^{n_1+1}\in[0;1]$ 
yields $\big(b_{n_1}^{(1)},\ldots,b_{n_2-1}^{(1)}\big)\in\THREE^{n_2-n_1}$
such that 
\[ \vec b^{(1)} \;:=\; \big(b_0^{(0)},\ldots,b_{n_1-1}^{(0)},
b_{n_1}^{(1)},\ldots,b_{n_2-1}^{(1)}\big)\in\THREE^{n_2} \]
has $|x_1-\sbin(\vec b^{(1)})|\leq2^{-n_2-1}$;
hence, in view of $|x_1-x_k|\leq2^{-n_2-1}$ according to Equation~\eqref{e:Multicont4},
$\vec b^{(1)}$ extends to signed binary expansions of any $x_k$, $k>1$. And so forth:
yielding a sequence of $\sigma$-names $\bar y^{(k)}\in\Cantor$ of $x_k\in\cball_{n_k+2}(x_{k-1})$ 
which agree on the first $n_k$ bits: $\bar y^{(k)}\in \cball_{n_{k}}(\bar y^{(k-1)})$.
\qed\end{enumerate}\end{proof}
%

\subsection{Representations over Compact Ultrametric Spaces}
\label{ss:Ultra}

Motivated by Example~\ref{x:Example},
Subsection~\ref{ss:Representations} had considered representations
over domains beyond Cantor space. And Theorem~\ref{t:Selection}
suggests which domains seem suitable for complexity considerations:
`almost' homogeneous compact ultrametric spaces.

The present section extends our notions of quantitative admissibility
(Section~\ref{s:Admissible2}) from representations over Cantor space 
to arbitrary compact ultrametric domains: based on quantitative continuity
of multifunctions according to Definition~\ref{d:Multicont}.
We then establish a quantitative \emph{Main Theorem} for multifunctions.

\medskip
Let us first generalize Definition~\ref{d:Admissible2}
from Cantor space with entropy $\id$ by Lemma~\ref{l:Entropy}a),
following the generic approach from Remark~\ref{r:Slack}.

\begin{definition}
\label{d:Admissible3}
Fix compact ultrametric space $(\calU,D)$ with entropy $\ENTROPY$.
\begin{enumerate}
\item[a)]
For (single-valued) representations $\xi:\subseteq \calU\twoheadrightarrow X$
and $\upsilon:\subseteq \calV\twoheadrightarrow Y$, 
a $(\xi,\upsilon)$-realizer of a \emph{multi}function $f:X\toto Y$ 
is a (single-valued) mapping $F:\subseteq\calU\to\calV$ 
satisfying any/all of the conditions from Observation~\ref{o:Restriction}c).
\item[b)]
Call representation $\xi:\subseteq \calU\twoheadrightarrow X$ \emph{polynomially admissible} ~iff~ 
it satisfies (b\,i) and (b\,ii), where:
\item[b\,i)]
$\xi$ has modulus of continuity $\mu$ such that $\ENTROPY\circ\mu \leq \Entropy\circ\poly$.
\item[b\,ii)]  
Every continuous surjective $\xi':\subseteq\calU\twoheadrightarrow X$ has $\xi'\reduceqP\xi$.
\item[c)]
Call representation $\xi:\subseteq \calU\twoheadrightarrow X$ \emph{polynomially standard} ~iff~ 
it satisfies (b\,i) $\ENTROPY\circ\mu \leq \Entropy\circ\poly$ and  
(b\,ii') $\xi'\reduceqP\xi$ for every (not necessarily surjective) $\xi':\subseteq\calU\to X$
and the following condition (c\,i):
\item[c\,i)]
$d(x,x')\leq2^{-\poly(n)}$ implies $d_\xi(x,x')\leq2^{-\mu(n)}$. 
\item[d)]
Call representation $\xi:\subseteq \calU\twoheadrightarrow X$ \emph{linearly admissible} ~iff~ 
it satisfies (d\,i) and (d\,ii), where:
\item[d\,i)]
$\xi$ has modulus of continuity $\mu$ with $\ENTROPY\circ\mu \leq \Entropy\circ\calO$.
\item[d\,ii)] 
Every continuous surjective $\xi':\subseteq\calU\twoheadrightarrow X$ has $\xi'\reduceqL\xi$.
\item[e)]
Call representation $\xi:\subseteq \calU\twoheadrightarrow X$ \emph{linearly standard} ~iff~ 
it satisfies (d\,i) $\ENTROPY\circ\mu \leq \Entropy\circ\calO$ and 
(d\,ii') $\xi'\reduceqL\xi$ for every (not necessarily surjective) $\xi':\subseteq\calU\to X$
and the following condition (e\,i):
\item[e\,i)]
$d(x,x')\leq2^{-\calO(n)} \quad\Rightarrow\quad d_\xi(x,x')\leq2^{-\mu(n)}$. 
\end{enumerate}
\end{definition}
Note that, compared to Definition~\ref{d:Admissible2},
$\Entropy$ has been replaced with $\mu$ in Items~(c\,i) and (e\,i):
which over Cantor space $\calU=\calC$ with $\ENTROPY=\id$ makes no difference because of (b\,i) and (d\,i), 
but for generic $\calU$ turns out as more convenient.
It does affect Items~(b) of the following generalized Theorem~\ref{t:Main1}, though:

\begin{theorem}
\label{t:Main2}
Fix compact ultrametric space $(\calU,D)$ with entropy $\ENTROPY$.
\begin{enumerate}
\item[a\,i)]
Suppose $(X,d)$ admits a polynomially standard representation $\xi:\subseteq\calU\twoheadrightarrow X$
and $\xi':\subseteq\calU\twoheadrightarrow X$ is polynomially admissible.
Then $\xi'$ is polynomially standard, too.
\item[a\,ii)]
Suppose $(X,d)$ admits a linearly standard representation $\xi:\subseteq\calU\twoheadrightarrow X$
and $\xi':\subseteq\calU\twoheadrightarrow X$ is linearly admissible.
Then $\xi'$ is linearly standard, too.
\item[b\,i)]
Fix a polynomially standard representation $\xi:\subseteq\calU\twoheadrightarrow X$ 
of metric space $(X,d)$ with entropy $\Entropy$. Consider $f:X\to Y$ for metric space $(Y,e)$.
If $f\circ\xi$ has modulus of continuity $\pmb{\mu}\circ\lambda$, then $f$ has modulus of continuity $\poly\circ\lambda$.
\item[b\,ii)]
If $\xi:\subseteq\calU\twoheadrightarrow X$ is linearly standard,
then $f$ has modulus of continuity $\calO\circ\lambda$.
\item[c)]
Let $\xi:\subseteq\calU\twoheadrightarrow X$ and $\upsilon:\subseteq\calU\twoheadrightarrow Y$
denote representations of metric spaces $(X,d)$ and $(Y,e)$ 
with respective entropies $\Entropy$ and $\theta$.
\item[c\,i)]
Suppose $\xi$ and $\upsilon$ are polynomially standard.
Consider $f:X\to Y$ with modulus of continuity $\lambda$.
Then $f$ admits a $(\xi,\upsilon)$-realizer $F:\dom(\xi)\to\dom(\upsilon)$
with modulus of continuity $\Lambda$ such that 
$\ENTROPY\circ\Lambda\circ\ENTROPY\circ\theta\leq \Entropy\circ\poly\circ\lambda\circ\poly$.
\item[c\,ii)]
Suppose $\xi$ and $\upsilon$ are linearly standard.
Consider $f:X\to Y$ with modulus of continuity $\lambda$.
Then $f$ admits a $(\xi,\upsilon)$-realizer $F$
with modulus of continuity $\Lambda$ such that 
$\ENTROPY\circ\Lambda\circ\ENTROPY\circ\theta\leq \Entropy\circ\calO\circ\lambda\circ\calO$.
\item[d)]
Let $\xi:\subseteq\calU\twoheadrightarrow X$ and $\upsilon:\subseteq\calU\twoheadrightarrow Y$
denote representations of metric spaces $(X,d)$ and $(Y,e)$ 
with respective entropies $\Entropy$ and $\theta$.
\item[d\,i)]
Suppose $\xi$ and $\upsilon$ are polynomially standard.
Let $f:X\to Y$ admit a $(\xi,\upsilon)$-realizer $F$ with modulus of continuity $\loinv{\ENTROPY}\circ\Entropy\circ\lambda$.
Then $f$ has modulus of continuity $\poly\circ\lambda\circ\loinv{\ENTROPY}\circ\theta\circ\poly$.
\item[d\,ii)]
Suppose $\xi$ and $\upsilon$ are linearly standard.
Let $f:X\to Y$ admit a $(\xi,\upsilon)$-realizer $F$ with modulus of continuity $\loinv{\ENTROPY}\circ\Entropy\circ\lambda$.
Then $f$ has modulus of continuity $\calO\circ\lambda\circ\loinv{\ENTROPY}\circ\theta\circ\calO$.
\end{enumerate}
\end{theorem}
Note that for Cantor space $\calU=\calC$ with $\ENTROPY=\id$
Theorem~\ref{t:Main2} recovers Theorem~\ref{t:Main1}.

We can now generalize Theorem~\ref{t:Standard}.
Theorem~\ref{t:Main2} supposes that the space(s) under consideration
be equipped with linearly/polynomially standard representations to begin with.
It remains to show that many spaces actually do admit such representations:

\begin{theorem}
\label{t:Standard2} 
Fix a compact \emph{ultra}metric space $(\calU,D)$ with minimum relative binary capacity $\CAPACITY=\CAPACITY_m(n)$.
Fix a compact metric space $(X,d)$ with binary entropy $\Entropy$,
and a strictly increasing positive integer sequence $n_k$, $k\in\IN$.
Inductively define $\mu(0):=\loinv{\CAPACITY}\circ\Entropy(1+n_0)$ and
\begin{equation}
\label{e:Standard2}
\mu(n)\;:=\;\loinv{\CAPACITY_{\mu(n_k)}}\circ\Entropy_{n_k}\big(1+n_{k+1}\big)\;\in\;\IN  
\quad \text{ for }n_{k}<n\leq n_{k+1} \enspace . 
\end{equation}
\begin{enumerate}
\item[i)]
Suppose that it holds 
\[ n_{K+1}\leq\calO(n_K) \quad\text{and}\quad \mu(n_K)\;\leq\; \loinv{\CAPACITY}\big(\Entropy(\calO(n_K))\big) \]
with respect to asymptotics $K\to\infty$ but ignoring constants depending on $(n_k)$.
Then there exists linearly admissible $\xi:\subseteq\calU\twoheadrightarrow X$ with modulus of continuity $\mu$.
\item[ii)] Suppose
$\displaystyle n_{K+1}\leq\poly(n_K) \quad\text{and}\quad
\mu(n_K)\;\leq\; \loinv{\CAPACITY}\big(\Entropy(\poly(n_K))\big)$ 
with respect to asymptotics $K\to\infty$ but ignoring constants depending on $(n_k)$.
Then there exists polynomially admissible $\xi:\subseteq\calU\twoheadrightarrow X$ with modulus of continuity $\mu$.
\end{enumerate}
\end{theorem}
%
Recall (Example~\ref{x:Entropy2}a) that Cantor space $\calU:=\calC$ has $\CAPACITY_m(n)=\max\{n-1-m,0\}$;
hence $\loinv{\CAPACITY_{m}}(N)=N+1+m$ and \[
\mu(n_{K+1}) \;=\; \Entropy_{n_K}\big(1+n_{K+1}\big)\:+\:1\:+\:\mu(n_K) \;=\cdots=\;
\Entropy\big(1+n_0)\:+\:\sum\nolimits_{k=0}^K \Entropy_{n_k}\big(1+n_{k+1}\big) \:+\: (K+1) \]
(almost) recovering Theorem~\ref{t:Standard}.

\medskip
The proof of Theorem~\ref{t:Standard2} is similar to that of Theorem~\ref{t:Standard}
on pages~\pageref{ss:Proofs1}ff, with the following changes: 
Recall the induction step $\xi_K\mapsto\xi_{K+1}$ with 
the mapping $\xi_K:\subseteq\TWO^{m_K}\to X$ satisfying conditions (I) and (II).
Defined inductively on Cantor space, 
$m_{K+1}:=m_{K}+\Entropy_{n_K}\big(1+n_{K+1}\big)$
guarantees $\dom(\xi_{K+1})$ to 
\begin{enumerate}
\item[III)] contain sufficiently many initial segments 
to encode the $\leq2^{m_{K+1}}$ centers of balls of radius $2^{-n_{K+1}-1}$ covering $\cball_{n_K}$ according to (I);
\item[IV)] encode in a discernible way in that these initial segments
(or rather: their extensions) have pairwise distances $>2^{-m_{K+1}}$
\end{enumerate}
and hence trivially satisfy the modulus condition~\eqref{e:Modulus},
while (II) makes sure that this property carries over to the induction step.

Generalizing Cantor space $\calC$ to $\calU$,
initial segments are replaced by points witnessing the capacity:
such that, in addition to unmodified (I) and (II),
\begin{enumerate}
\item[III')] $\dom(\xi_{K+1})\subseteq\calU$ has cardinality $\geq2^{m_{K+1}}$ 
\item[IV')] and any two points in $\dom(\xi_{K+1})$ have pairwise distances $>2^{-\mu(n_{K+1})}$
\end{enumerate}
and hence again trivially satisfy the modulus condition~\eqref{e:Modulus}
with $\mu$ according to Equation~\eqref{e:Standard2}:
Note that arguments in $\cball_{m_K}(x)\subseteq\calU$ having pairwise distance $>2^{-\mu(n_{K+1})}$ from each other
also have such distance to arguments in $\cball_{\mu(n_K)}(x')$ for $D(x,x')>2^{-\mu(n_K)}$ due to $(\calU,D)$ being ultrametric.
\qed

\medskip\noindent
Proceeding from Cantor to other ultrametric spaces 
provides a promising approach to the open problem
of defining (!) computational complexity for higher types \cite{Royer}:

\begin{example}
\label{x:Higher}
Returning to Example~\ref{x:Example},
consider the application functional
\[ \App \;:\; \LIP([0;1],[0;1])\times[0;1] \;\ni\; (f,x) \;\mapsto \; f(x) \;\in\; [0;1] \]
which is non-expansive \ie 1-Lipschitz.
According to Lemma~\ref{l:Entropy}e+f), 
no Cantor space representation of $\dual{[0;1]}=\LIP([0;1],[0;1])$
can have a sub-exponential modulus of continuity.

In fact, for every $\Delta:\subseteq\Cantor\twoheadrightarrow \dual{[0;1]}=\LIP([0;1],[0;1])$
and with the Cartesian product representation $\Delta\times\delta$ of $\dual{[0;1]}\times[0;1]$,
no $(\Delta\times\delta,\delta)$-realizer $F:\subseteq\Cantor\to\Cantor$
can have a sub-exponential modulus of continuity:
the seemingly simple application functional is not computable in polynomial time---when
using Cantor space representations.

The representation $\deltabox$ from \cite[\S3.4]{KC12} on the other hand
is (or more precisely: can be rewritten and restricted to be) defined 
over the compact ultrametric space $\calU:=\dual{(\TWO^*)}$, 
which according to Example~\ref{x:Entropy2}b)
has exponential binary entropy $\ENTROPY(n)=2^n$ and 
relative binary capacity $\CAPACITY_m(n)=2^{n-1}-2^m$ for $n>m$.
Hence, over $\dual{(\TWO^*)}$ instead of $\Cantor$, 
Lemma~\ref{l:Entropy}e+f) does not prohibit a representation
of $\dual{[0;1]}$ from having polynomial modulus of continuity.
And indeed, $\App$ is $(\deltabox\times\delta,\delta)$-computable
in polynomial time \cite[Lemma~4.9]{KC12}.
\end{example}


\section{Concluding Comments}
\label{s:Conclusion}

The present work has boiled down the original Question~\ref{q:Main}
to Definition~\ref{d:Admissible2}, thoroughly justified 
with the many benefits of this new notion of quantitative admissibility:
existence (Theorem~\ref{t:Standard}) and closure properties (Theorem~\ref{t:Cartesian})
and quantitative strengthenings (Theorems~\ref{t:Main1}) 
of the classical qualitative \emph{Main Theorem} (Fact~\ref{f:KreitzWeihrauch})
that generalize and systematize (Theorem~\ref{t:Main2}) previous results about known ad-hoc encodings (Corollary~\ref{c:Main1}).

\medskip
It remains to further investigate the following purely mathematical question:

\begin{question}
What are sufficiently `rich' classes/categories of compact metric spaces $(X,d)$ with relative binary entropy $\Entropy=\Entropy_m(n)$
which admit linearly/polynomially admissible/standard representations
in the sense of Definition~\ref{d:Admissible2}?

What are sufficiently `rich' classes/categories of compact metric spaces 
that satisfy the hypothesis of Theorem~\ref{t:Standard}, \ie,
for which there exist strictly increasing positive integer sequences $n_k$
such that $\mu(n_K):=\Entropy(1+n_0)+\sum_{k=0}^K \Entropy_{n_{k}}\big(1+n_{k+1}\big)$ 
satisfies 
\begin{enumerate}
\item[i)] $n_{k+1}\leq\calO(n_k)$ and 
\item[ii)] $\mu(n_K)\leq \Entropy_X\big(\calO(n_K)\big)$ ?
\end{enumerate}
\end{question}
\medskip\noindent
Throughout this work we have focused on representing compact metric spaces $(X,d)$.
These admit complexity analyses in dependence on one integer parameter $n$ only; cmp. \cite{Sch04}.
Generalizing to sigma-compact metric spaces $X=\bigcup_k X_k$ with compact $X_k$
leads to \emph{parameterized} complexity, with additional dependence on secondary parameter $k\in\IN$
\cite{EikeFlorian}. 

\addcontentsline{toc}{section}{\refname}
\bibliographystyle{alpha}
\bibliography{cca,multicont}

\end{document}